\documentclass[leqno,12pt]{article}
\usepackage{amssymb,amsmath,amscd}

\makeatletter
\@addtoreset{equation}{subsection} 

\renewcommand{\@begintheorem}[2]{
\begin{trivlist}\it \item[\hspace{\labelsep}{\bf #1\ #2\ }]}

\renewcommand{\@opargbegintheorem}[3]{
\begin{trivlist}\it \item[\hspace{\labelsep}{\bf #1\ #2\ (#3)\ }]}

\renewcommand{\@endtheorem}{\end{trivlist}}
\makeatother

\newtheorem{teo}[equation]{Theorem}
\newtheorem{cla}[equation]{Claim}
\newtheorem{lem}[equation]{Lemma}

\newtheorem{pro}[equation]{Proposition}
\newtheorem{cor}[equation]{Corollary}

\newtheorem{teo1}{Theorem}

\newenvironment{ea*}{\begin{eqnarray*}}{\end{eqnarray*}}

\newenvironment{rem}{ \refstepcounter{equation}\par\noindent
 {{\bf Remark\/}  \bf \thesubsection.\arabic{equation}\ \ }}{\par}
 
\newenvironment{dfn}{ \refstepcounter{equation}\par\noindent
 {{\bf Definition\/} \bf \thesubsection.\arabic{equation}\ \ }}{\par}

\begin{document}

    
\def\PP{{\mathcal{P}}}
\def\RR{{\mathcal{R}}}
\def\AA{{\mathcal{A}}}
\def\XS{{  \Xi_{sing}   }}
\def\TS{{  \Theta_{sing}   }}
\def\LL{{\mathcal{L}}}
\def\RR{{\mathcal{R}}}
\def\UU{{\mathfrak{U}}}
\def\CT{{  \widetilde{C}  }}
\def\NT{{  \widetilde{N}  }}
\def\DT{{  \widetilde{D}  }}
\def\THT{{ \widetilde{\Theta} }}
\def\PR{{  \PP (\tilde{C},C,\Xi)  }}
\def\AT{{  (A, \Theta)  }}
\def\OO{{\mathcal{O} }}
\def\XE{{  \Xi_{ex}  }}
\def\XST{{ \Xi_{st}  }}
\def\SS{{  \mathbb{C}[t]/(t^{n+1})  }}
\def\lra{{  \longrightarrow  }}
\def\ra{{  \rightarrow  }}
\def\PROJ{{  \mathbb{P}     }}
\def\PR{{    \mathbb{P}^{1} }}
\def\PRR{{   \mathbb{P}^{2} }}
\def\PRRR{{  \mathbb{P}^{3} }}
\def\PRRRR{{ \mathbb{P}^{4} }}

\def\BB{{\mathcal{B}  }}
\def\CC{{ \mathcal{C} }}
\def\MM{{ \mathcal{M} }}
\def\TT{{ \mathcal{T} }}
\def\MT{{ \widetilde{M} }}
\def\MH{{ \widehat{M} }}

\def\aa{{ \alpha }}
\def\ss{{ \sigma }}
\def\tt{{ \tau }}
\def\ts{{ \tilde{s} }}
\def\tm{{ \tilde{m} }}
\def\hm{{ \hat{m} }}
\def\ta{{ \tilde{a} }}
\def\ll{{ \lambda }}
\def\mm{{ \mu     }}  
\def\ffi{{\varphi }}
\def\pp{{ \partial }}
\def\FF{{ \mathcal{F} }}
\def\GG{{ \mathcal{G} }}
\def\CX{{ \mathbb{C} }}
\def\ZZ{{ \mathbb{Z} }}
\def\RE{{ \mathbb{R} }}
\def\QQ{{ \mathbb{Q} }}
\def\MA{{ \textnormal{\textbf{M}} }}

\def\CH{{ \check{\textnormal{C}} }}
\def\SCH{{ \check{\mathcal{C}} }}
\def\mam{{\mathfrak{m}}}

\def\KER{{ \textnormal{ker} }}
\def\IM{{ \textnormal{im} }}
\def\PIC{{\textnormal{Pic}  }}
\def\DIV{{\textnormal{Div}  }}
\def\SPEC{{\textnormal{Spec }  }}
\def\MUL{{ \textnormal{mult} }}
\def\COK{{ \textnormal{coker} }}
\def\RES{{ \textnormal{Res} }}
\def\DET{{ \textnormal{det} }}
\def\QED{{ \hspace*{\stretch{1}}$\square$  }}
\def\PRF{{\noindent\emph{Proof}}}
\def\remskip{{ \vskip .12in }}
\def\RK{{ \textnormal{rank} }}
\def\sing{{\textnormal{Sing }  }}
\def\singd{{\textnormal{Sing}  }}
\def\supp{{\textnormal{supp}  }}
\def\cliff{{\textnormal{Cliff}  }}


\title{Singularities of the Prym Theta Divisor}
\author{Sebastian Casalaina-Martin\footnote{The author 
was partially supported by a VIGRE fellowship from NSF grant DMS-98-10750}}
\maketitle

\begin{quote}
\begin{center}\textbf{Abstract} \end{center}

For the Jacobian of a curve, the Riemann singularity theorem gives a
geometric interpretation of the singularities of the theta divisor in
terms of special linear series on the curve.  This paper proves an
analogous theorem for Prym varieties.  Applications of this theorem to
cubic threefolds, and Prym varieties of dimension five, are also
considered. \end{quote}

\section*{Introduction}

A principally polarized abelian variety (ppav) can be studied
through the geometry of its theta divisor. While in general this
geometry is not well understood, one can simplify the problem by
focusing on ppav's related to curves.  This paper will consider such varieties defined over the complex numbers.  Among the most studied
examples are Jacobians: given a smooth curve $C$ of genus $g$, the
Jacobian of $C$ is the ppav of dimension $g$ defined as $JC=H^0(C,
\omega_C)^*/H_1(C,\ZZ)$.  The polarization is given by a theta
divisor $\Theta$, whose geometry is closely related to the curve
$C$.  The Abel-Jacobi theorem allows us to identify $JC$ with
$\PIC^{g-1}(C)$, the space of isomorphism classes of line bundles
over $C$ of degree $g-1$, and it is often convenient to make this
identification when studying $\Theta$.  A fundamental result is
Riemann's singularity theorem, which states that $\MUL_x
\Theta=h^0(C,L_x)$, where $L_x$ is a line bundle of degree $g-1$
in the isomorphism class associated to $x$. 

There is a close connection between the singular points of the
theta divisor and the canonical image of $C$:  the tangent cone to
a general singular point of $\Theta$ contains the canonical image
of $C$, and a theorem of Green's \cite{green} implies that for a smooth curve
$C$ of genus $g\ge 4$, with no $g^1_2,\ g^1_3$ or $g^2_5$, the
canonical image of $C$ is cut out by the tangent cones to double
points of $\Theta$.  In general, a theorem due to Torelli states
that the Jacobian $(JC, \Theta)$ uniquely determines the curve
$C$. 

Another commonly studied ppav is the Prym variety of a connected
\'etale double cover of a smooth curve.  The study of these
varieties goes back to Riemann, and the geometry of the Prym theta
divisor is in many ways parallel to that of the Jacobian theta
divisor.  Recall that associated to $\pi:\CT \ra C$, a connected
\'etale double cover of a smooth curve $C$ of genus $g$, is an
involution $\tt:\CT \ra \CT$, which induces an involution on
$H^0(\CT,\omega_\CT)^*$ and $H_1(\CT,\ZZ)$. 
Denoting by $(H^0(\CT,\omega_\CT)^*)^-$ and $H_1(\CT,\ZZ)^-$
the negative eigen spaces of the involution, 
the Prym variety $P$
associated to such a cover is the $g-1$ dimensional abelian
subvariety of $J\CT$ defined as 
$P=(H^0(\CT,\omega_\CT)^*)^-/H_1(\CT,\ZZ)^-$.  As in the case of
Jacobians, it is convenient when studying Pryms to identify $J\CT$
with $\PIC^{2g-2}(\CT)$. In this case, $P$ can be described set
theoretically as $$P= \{ L \in \PIC^{2g-2}(\CT) \ | \
\textnormal{Norm}(L)=\omega_C, \ h^0(L) \equiv 0 \
(\textnormal{mod } 2) \}.$$ There is a principal polarization on
$P$ given by a theta divisor $\Xi$;  as a set $\Xi=\{ L \in P \ |
\ h^0(L) \ge 2 \}$.  In this paper, we will prove an analog of the
Riemann singularity theorem in the case of Pryms; i.e. we will
relate the multiplicity of a point $x\in \sing \Xi$, to the
dimension of special linear systems on $\CT$ and $C$. 

The study of this question goes back to Mumford, who proved in
\cite{mum1} that if $x\in \Xi$, then $\MUL_x \Xi \ge h^0(L_x)/2$,
and that if in addition $h^0(L_x)=2$, then $x \in \sing \Xi$ if
and only if $L_x=\pi^*(M)\otimes \OO_{\CT}(B)$, where $M$ is a
line bundle on $C$ such that $h^0(M)=2$, and $B$ is an effective
divisor on $\CT$. The proof of the second statement is based on
Kempf's generalization of the Riemann singularity theorem
\cite{kempf}.  Smith and Varley used a similar method in
\cite{sv2} to prove the following theorem extending Mumford's
assertions:

 \begin{teo1}[Smith-Varley \cite{sv2}]\label{teosv} Let $\pi: \CT
\rightarrow C$ be a connected \'etale double cover of a smooth
curve $C$ of genus $g$, and let $(P,\Xi)$ be the associated Prym
variety.  If $x \in \Xi$ corresponds to a line bundle $L\in
\PIC^{2g-2}(\CT)$, and $C_x\THT$ is the tangent cone to $\THT$ at
$x$, then the following are equivalent: \begin{itemize}

\item[\textnormal{(a)}] $T_x P \subseteq C_x \THT$;

\item[\textnormal{(b)}] $\MUL_x \Xi >h^0(L)/2$;

\item[\textnormal{(c)}] $L=\pi^*M \otimes \OO_{\CT}(B)$,
$h^0(C,M)>h^0(\CT,L)/2$, $B \ge0$, and $B \cap \tau^*B =
\emptyset$. \end{itemize} 
Furthermore, if \textnormal{(c)} holds,
then $M$ and $\OO_{\CT}(B)$ are unique up to isomorphism, and
$h^0(\CT,\OO_{\CT}(B))=1$. 

\end{teo1}

In this paper, we will complete the analog of the Riemann
singularity theorem by determining the exact multiplicity of
singular points of this type: 

 \begin{teo1}\label{teo2} 
In the above notation, suppose that $x$
is a singular point of $\Xi$, corresponding to a line bundle $L\in
\PIC^{2g-2}(\CT)$ such that $L=\pi^*M \otimes \OO_{\CT}(B)$, with
$h^0(C,M) > h^0(\CT,L)/2$, $B \ge0$, and $B \cap \tau^*B =
\emptyset$.  Then $\MUL_x \Xi=h^0(M)$. 

\end{teo1}

To prove the theorem, we will consider a deformation of the line
bundle $L$, and then relate $\MUL_x \Xi$ to the obstruction to
lifting sections of $L$ to sections of the deformation.  In the
course of proving the theorem, we will give a short proof of Smith
and Varley's theorem.

\remskip 
More generally, (see \cite{birlange}, p.371), a ppav
$(Z,\Xi)$ is called a Prym-Tjurin variety or generalized Prym
variety if there is a smooth curve $C$ with Jacobian $(JC,\Theta)$
such that $Z$ is an abelian subvariety of $JC$ with $\Theta \cap
Z= e \cdot \Xi$.  Necessarily, $e$ is the exponent of $Z$ in $JC$. 
Every ppav of dimension $g$ is a
Prym-Tjurin variety of exponent $3^{g-1}(g-1)!$, and it is
possible that there are situations other than the exponent $2$
case examined in this paper, where these techniques can be used to
compute the multiplicity of a singular point of $\Xi$.

\remskip
For a ppav $(A,\Theta)$, 
let $\singd_k\Theta=\{x\in \sing\Theta\ |\ \MUL_x \Theta\ge k\}$. 
A result of Koll\'ar's \cite{kbound} shows that if
$\dim(A)=d$, then $\dim(\singd_k\Theta)\le
d-k$.  Generalizing a result of Smith and Varley \cite{svbound}, Ein and Lazarsfeld \cite{elbound} showed  that 
$\dim(\singd_k\Theta)=d-k$ only if $(A,\Theta)$ splits as a 
$k$-fold product.  In particular, for an irreducible ppav of dimension $d$, and a point $x\in \Theta$, $\MUL_x\Theta\le d-1$.  
For the Jacobian of a smooth curve of genus $g$, applying the Riemann singularity theorem, and Martens' theorem 
\cite{martens}, one can see that these bounds are not optimal; in fact
$\dim(\singd_k\Theta)\le g-2k+1$, with equality holding only if $C$ is hyperelliptic.  This implies in particular that 
$\MUL_x\Theta\le (g+1)/2$.

For a Prym variety  associated to a connected \'etale double cover of a smooth curve $C$ of genus $g$, and 
$x\in \sing \Xi$, the bounds given above yield that
 $\MUL_x \Xi\le g-2=\dim(P)-1$.  
Although in \cite{mum1} Mumford used a strengthened version of Martens' theorem to prove statements about the dimension of the singular locus of the Prym theta divisor,
without the results of 
Theorem \ref{teo2}, no bound could be given on the multiplicity of these points.  The implications of Theorem
\ref{teo2} for the singular locus of the Prym theta divisor are described in Corollary \ref{upbound}.  In particular, 
it is shown that for an irreducible Prym variety 
$(P,\Xi)$, and a point $x\in \Xi$, $\MUL_x\Xi \le 
(\dim(P)+1)/2$.

\remskip

The rich connection between singularities of the Jacobian theta
divisor and the canonical image of the curve, is reflected in the
Prym case by the connection between singularities of the Prym
theta divisor, and the Prym canonical image of the base curve $C$.  A result of Tjurin (\cite{t},  Lemma 2.3, p.963), generalized by
Smith and Varley in (\cite{sv1}, Proposition 5.1), shows that the
Prym canonical image of $C$ is contained in the tangent cone to
$\Xi$ at $x$, for all singular double points $x$ such that $T_x P
\nsubseteq C_x \THT$.  In addition, there are many such points for
curves of high genus.  Analogous to Green's theorem for Jacobians \cite{green},
a primary open question for Prym varieties is to determine when
the quadric tangent cones to $\Xi$ cut out the Prym canonical
image of $C$.  In \cite{debarre}, Debarre showed that this is true
for general curves of genus $\ge 8$. In other words, the Prym map
$ \PP: \RR_{g}\ra \AA_{g-1}$, taking a connected \'etale double
cover of a smooth curve to its associated Prym variety, is
generically injective.  Friedman and Smith have shown in \cite{fs}
that the Prym map is generically injective for $g\ge 7$. On the
other hand, unlike the case of Jacobians, there are examples in
every dimension of Pryms arising from non-isomorphic double covers
(see \cite{sv3}, p.237), and the question remains open to
determine exactly which Pryms arise from a unique double cover of
curves.

As a step towards understanding this question, we will consider
the intersection of the $k$-secant variety of the Prym canonical image
of $C$, that is the variety of $k$-dimensional secants
(c.f. Section \ref{secsec}), with the tangent cone to $\Xi$ at singular points.  In
particular, if $\CT$ is not hyperelliptic, $\Psi$ is the Prym
canonical morphism of $C$, and $C_x \Xi$ is the projective tangent
cone to $\Xi$ at $x$, then we will show the following: suppose
$x\in \sing \Xi$ corresponds to the line bundle $L\in
\PIC^{2g-2}(\CT)$, and $h^0(L)=2n$, then the $(n-1)$-secant
variety to $\Psi(C)$ is not contained in $C_x \Xi$, while the
$(n-2)$-secant variety is contained in $C_x\Xi$.  As a
consequence, $\Psi(C)\subseteq C_x \Xi$ if and only if $h^0(L)\ge
4$.  Thus, if one hopes to recover the curve $C$ as the base locus
of quadric tangent cones to $\Xi$, one must exclude the tangent
cones at points with $h^0(L)=2$.  This was suspected to be true
since Smith and Varley \cite{sv3} have observed that 
if $(P,\Xi)$ is the Jacobian of a non hyperelliptic curve
and $x$ is a generic double point such that $h^0(L)=2$, then $\Psi(C)\nsubseteq
C_x\Xi$ (c.f. Remark \ref{remtjurin}).

\remskip

Prym varieties also arise in the study of conic bundles, and in
particular, in the study of cubic threefolds. Mumford stated in
\cite{mum1} that the intermediate Jacobian of a smooth cubic
threefold in $\PRRRR$ is isomorphic to the Prym variety of a
connected \'etale double cover of a smooth plane quintic.  Using
this description of the intermediate Jacobian, he stated the
following theorem: if $X$ is a smooth cubic threefold in $\PRRRR$
with intermediate Jacobian $(JX,\Theta)$, then $\sing
\Theta=\{x\}$, $\MUL_x \Theta=3$, and moreover, $C_x \Theta \cong
X$.  It follows from this that $JX$ determines $X$ up to
isomorphism, and $X$ is irrational; both statements were first
proven by Clemens and Griffiths in \cite{cg}.

In \cite{cmf}, a converse to Mumford's theorem was proven:  if
$(A,\Theta)$ is a ppav of dimension $5$, $\sing \Theta=\{x\}$ and
$\MUL_x\Xi =3$, then $(A,\Theta)$ is isomorphic to the
intermediate Jacobian of a smooth cubic threefold.  If one removes
the condition that $\sing \Theta =\{x\}$, and requires instead
the weaker condition that exactly one of the singular points of $\Theta$ has multiplicity $3$,
then it was shown that the only other possibility is that
$(A,\Theta)$ is isomorphic to $JC$ or $JC\times JC'$ for some
hyperelliptic curves $C$ and $C'$ (i.e. curves having a line
bundle $L$ such that $\deg(L)=h^0(L)=2$). One would like to have a
complete description of all ppav's of dimension five whose theta
divisor has a triple point.  The following theorem is a
consequence of Theorem \ref{teo2}: 

\begin{teo1}\label{teogen6}
Let $(P,\Xi)$ be a Prym variety of dimension five, and let $x \in
\sing \Xi$.  If $\MUL_x \Xi =3$, then $(P,\Xi)$ is a hyperelliptic
Jacobian, the product of two hyperelliptic Jacobians, or the
intermediate Jacobian of a smooth cubic threefold. \end{teo1}

The proof relies on a theorem of Mumford's \cite{mum1} concerning
Prym varieties of hyperelliptic curves, and on a theorem due to
Beauville \cite{b2} (see also Donagi and Smith \cite{ds}),
concerning Prym varieties of plane quintics. \remskip

The Prym map $\PP : \RR_6 \ra \AA_5$ has dense image, and so to
extend Theorem \ref{teogen6} to all ppav's of dimension five, it
remains only to check the statement on the boundary of the image. 
Beauville has shown in \cite{b1} that the ppav's on the boundary
correspond to Prym varieties of admissible double covers of stable
curves of genus $6$.  It is reasonable to expect that Theorem
\ref{teo2} can be extended to this case, at least for curves of
low genus.  Consequently, one should be able to describe all ppav's
of dimension five whose theta divisors have triple
points.  This is work in progress. 

\remskip

The outline of the paper is as follows.  Section 1 focuses on how
to calculate the multiplicity of points of $\sing \Xi$.  Section 2
concerns some general results that will be useful for computations
in later sections. In Section 3 we give a short proof of Smith and
Varley's theorem. In Section 4 we prove Theorem \ref{teo2}, and in
Section 5 we prove some immediate consequences, including Theorem
\ref{teogen6}. Section 6 establishes the connection between the
computations made in Sections 4 and 5, and the Prym canonical
image of the base curve. In particular, we examine the secant
variety to the Prym canonical curve.  We also give a brief
description of the equation defining the tangent cone to $\Xi$ at
certain singular points.

\remskip

The results in this paper appear as part of the author's Ph.D.
thesis. I would like to take this opportunity to thank my advisor
Robert Friedman for his generous help. I would like to thank
Roy Smith and Robert Varley for their detailed comments, and for providing
me with copies of their unpublished work; this was particularly
helpful in formulating the statement of Theorem \ref{teo2}. 
I would also like to thank the referee for suggesting many improvements, especially the strengthening of 
Corollary \ref{upbound} to include statements about the dimension of the singular locus of the Prym theta divisor.

\section{Theta Divisors}

In this section we state the key results from \cite{cmf}, which we
will need in what follows.  The proofs of these facts will be
omitted, except in the cases where certain generalizations are
needed. 

\subsection{Preliminaries on theta divisors}

Let $S$ be a scheme, $C$ be a smooth, connected, complete curve, and $\LL$ be a line
bundle over $C \times S$, of relative degree $g-1$.  Then
a result of Grothendieck (\cite{EGA} 6,7) as formulated by Mumford (\cite{mumab} Theorem p.46, or \cite{hart}
III Lemma 12.3) gives, locally on $S$, a complex of locally free
$\OO_{S}$-modules of the same rank, $d:\CC^{0} \ra \CC^{1}$, whose
cohomology is $R^{0}\pi_{2*}\LL$ in dimension zero, and
$R^{1}\pi_{2*}\LL$ in dimension one.  If $(\det d)$ is not a zero
divisor, then $(\det d)$ is an effective Cartier divisor
which is independent of the choice of the complex $\CC^{\bullet}$,
and hence defines a global effective divisor $\Theta_S$ on $S$,
which satisfies the following:

\begin{teo}\label{theta}
In the above notation, $\Theta_{S}$ is an effective nonzero
Cartier divisor on $S$, and satisfies: \begin{itemize}
    
\item[\textnormal{(a)}] 
the support of $\Theta_{S}$ is equal to the set of $s\in S$ 
such that $h^{0}(C;\LL_{s}) \ne 0$;

\item[\textnormal{(b)}]  if $S= \PIC^{g-1}(C)$ and $\LL$ is a Poincar\'e line 
bundle, then $\Theta_{S}$ is the usual theta divisor;

\item[\textnormal{(c)}] the construction is functorial:  if $f:S' \ra S$ is a 
morphism, and $\LL'=( \textnormal{Id} \times f)^{*}\LL$, then
$\Theta_{S'}=f^{*}\Theta_{S}$;

\item[\textnormal{(d)}] if $S$ is smooth and dim $S=1$, then 
$\Theta_{S}=
\sum_{s\in S}\ell((R^{1}\pi_{2*}\LL)_{s})\cdot s$,
where $(R^{1}\pi_{2*}\LL)_{s}$ refers to the stalk at $s$.
\end{itemize}
\end{teo}
\PRF.  (a)-(c) are standard.  A reference for the proof of (d) is Friedman and Morgan \cite{fm} Proposition 3.9, p.384.\QED
\remskip

Let us now restrict to the case that $S$ is a smooth curve with
$s_0 \in S$. Let $t$ be a local coordinate for $S$ centered at
$s_{0}$ which only vanishes there, and set $S_k = \SPEC
\CX[t]/t^{k+1}$.  For each $k$, there is a map $S_k \rightarrow
S$, so that if we set $C_{k}= C \times S_k$, then there are
induced maps $C_k \rightarrow C \times S$. For example, $C_0=C$,
and $C_0\rightarrow C\times S $ is the inclusion of the fiber over
$s_0$. Finally, let $\LL_{k}$ be the restriction of $\LL$ to
$C_{k}$.  It follows that $\LL_0=L$ is the restriction of $\LL$ to
$C \times \{ s_{0} \}$, and $\LL_k= \LL/ t^{k+1}\LL$. 

\begin{lem}[\cite{cmf}, 1.5]\label{lemR1}
For all $k$, $\ell(H^{0}(C_{k}, \LL_{k})) \le \ell(H^{0}(C_{k+1},
\LL_{k+1}))$. Furthermore, there is an $N\in \ZZ$ such that for
all $k\ge N$, $\ell(H^{0}(C_{k}, \LL_{k}))$ is independent of $k$
and $$\ell(H^{0}(C_{k}, \LL_{k}))=\ell((R^{1}\pi_{2*}\LL)_{s_{0}})
=\MUL_{s_0}\Theta_S.$$ \QED 
\end{lem}

We can be more explicit about the value of $N$.  There is an exact 
sequence 
\begin{equation}\label{sesk}
0 \lra t\LL_{k} \lra \LL_{k} \lra L \lra 0,
\end{equation}
where $t\LL_{k}\cong \LL_{k-1}$, and the obvious surjection 
$\LL_{k} \ra \LL_{k-1}$ induces a commutative diagram

\begin{equation}\label{diagN}
\begin{CD}
0 @>>> H^{0}(\LL_{k-1}) @>>> H^{0}(\LL_{k})   @>>> H^{0}(L) @>\partial_k>>      H^{1}(\LL_{k-1}) \ldots \\
@.     @VVV                  @VVV                  @|                           @VVV \\ 
0 @>>> H^{0}(\LL_{k-2}) @>>> H^{0}(\LL_{k-1}) @>>> H^{0}(L) @>\partial_{k-1}>>  H^{1}(\LL_{k-2})\ldots 
\end{CD}
\end{equation}

\remskip
\begin{lem}[\cite{cmf}, 1.6]\label{inj}
Suppose in the above notation that $\partial_{N+1}$ is injective
for some $N$.  Then for all $k \ge N$, the natural inclusion
$t^{k-N}\LL_{k} \subseteq \LL_{k}$ induces an equality
$H^{0}(t^{k-N}\LL_{k})=H^{0}(\LL_{N})$.  In particular
$\ell(H^{0}(\LL_{k}))= \ell(H^{0}(\LL_{N})) $ for all $k \ge N$. 
\QED \end{lem}

One would like to have a way of computing $\ell(H^{0}(\LL_{k}))$. 
Define $W_k$ to be the image of the map $H^0(\LL_k) \ra H^0(L)$
induced by the exact sequence (\ref{sesk}), and let
$d_k=\dim(W_k)$.  We will say that a section $s\in H^0(L)$ lifts
to order $k$ if $s\in W_k$.  It is clear from the commutivity of
the diagram (\ref{diagN}) that for all $k$, $W_{k+1}\subseteq W_k$
and hence $d_{k+1}\le d_k$. 

\begin{lem}\label{lemlength}
In the notation above, 
$\ell(H^{0}(\LL_{k}))= 
        \sum_{i=0}^{k}d_{i}$.
\end{lem}

\PRF.  This follows by induction on $k$ using the following exact
sequence \[0 \lra H^0(\LL_i)\lra H^0(\LL_{i+1})\lra H^0(L)
\stackrel{\partial_{i+1}}{\lra} H^1(\LL_i)\lra \dots \] Indeed,
$\ell(H^0(\LL_{i+1}))=\ell(H^0(\LL_{i}))+d_{i+1}$. \QED

\subsection{Obstructions to lifting}

For an appropriate affine cover $\{ U_i \}$ of $C$, we may assume
that $L$ has transition functions $\ll_{ij}$ and that the
transition functions for $\LL$ are of the form
$\ll_{ij}(t)=\ll_{ij}(1+\sum_{k=1}^\infty \aa_{ij}^{(k)}t^k)$.  By
definition, these satisfy the condition $\ll_{ik}(t) = \ll_{ij}(t)
\ll_{jk}(t)$, and it follows that the cochain $\xi= \aa_{ij}^{(1)}$ 
is a cocycle in $H^1(\OO_C)$.  Likewise set

\[ \ll_{ij;N}(t)=\ll_{ij}(1+\sum_{k=1}^N \aa_{ij}^{(k)}t^k). \]
Assume that $s\in H^0(L)$ and that $s_{N-1}$ is a lifting of $s$
to a section of $\LL_{N-1}$.  Then using the trivialization over
the open cover $\{ U_i \}$ we have \[ s_{i;N-1}= \sum_{k=0}^{N-1}
\ss_{i}^{(k)}t^k \] for some functions $\ss_i^{(k)}\in
\OO_C(U_i)$, with $s_{i;N-1}= \ll_{ij;N-1}(t)s_{j;N-1}$ on $(U_i
\cap U_j) \times \SPEC \CX[t]/(t^N)$. The section $s_{N-1}$ lifts
to a section $s_N$ if and only if there exists $\ss_i^{(N)} \in
\OO_C(U_i)$ such that, if we set $ s_{i;N}=\sum_{k=0}^{N}
\ss_{i}^{(k)}t^k $ then $s_{i;N}= \ll_{ij;N}(t)s_{j;N}$ on $(U_i
\cap U_j) \times \SPEC \CX[t]/(t^{N+1})$.  Since $s_{i;N-1}$ is
already a section of $\LL_{N-1}$, this is equivalent to the
condition

\[ \ss_i^{(N)} = \ll_{ij}\ss_j^{(N)}+\sum_{k=0}^{N-1}\ll_{ij}\aa_{ij}^{(N-k)}
\ss_j^{(k)}.\]

Let $\gamma_{N}(s_{N-1})$ be the $1$-cochain defined by
$\sum_{k=0}^{N-1}\ll_{ij}\aa_{ij}^{(N-k)}\ss_j^{(k)}$; i.e. the
obstruction to lifting $s_{N-1}$ to order $N$.  I claim
$\gamma_{N}(s_{N-1})$ is a $1$-cocyle in $H^1(L)$. Indeed, let
$\partial_k^{\circ}$ denote the map $H^0(\LL_{k-1})\ra$ $ H^1(L)$
induced from the exact sequence \begin{equation} 0 \lra L \lra
\LL_k \lra \LL_{k-1} \lra 0. \end{equation} A computation in the
$\CH$ech complex will then show: 

\begin{lem} Suppose $s_{N-1} \in H^0(\LL_{N-1})$.  Then
$\gamma_N(s_{N-1})=\partial_{N-1}^{\circ}(s_{N-1}) $ $\in H^1(L)$, and
thus $s_{N-1}$ lifts to a section $s_N \in H^0(\LL_N)$ if and only
if $\gamma_N(s_{N-1})$ $=0$ in $H^1(L)$.\QED 
\end{lem}

Computing these obstructions is the central step in the proofs of
the main theorems.  In these proofs, we will be restricting our
attention to a particular class of deformations described in the
next section, and in that case we will write down explicit
formulas for the first and second order obstructions.  We will
also outline a particular technique for determining their class in
$H^1(L)$. The basic idea is illustrated by the following lemmas
regarding first order lifts. 

\begin{lem}
Let $\xi\in H^1(\OO_C)=\textnormal{Ext}^1(L,L)$ be the extension
class corresponding to $\LL_1$.  Then $\gamma_1(s)=s\cup \xi \in
H^1(L)$, where the cup product is $H^0(L)\otimes H^1(\OO_C)\ra
H^1(L)$.\QED
\end{lem}

\begin{lem}
Let $D$ be an effective divisor on $C$ and let $\partial$ be the
coboundary map $H^0(\OO_D(D))\ra H^1(\OO_C)$ induced by the short
exact sequence $$ 0\ra \OO_C\ra \OO_C(D)\ra \OO_D(D)\ra 0. $$
Suppose that $\xi \in H^1(\OO_C)$ is of the form $\partial(t)$ for
some $t\in H^0(\OO_D(D))$.  Then $s\cup \xi=\partial_L(s\cdot t)$,
where $s\cdot t$ is the section of $L(D)|_D$ given by taking the
cup product of $s$ and $t$, and $\partial_L$ is the coboundary
homomorphism arising from $$ 0\ra L\ra L(D)\ra L(D)\otimes
\OO_D\ra 0.  $$ \QED
\end{lem}

We now consider the following useful observation. Suppose that $p$
is a point of $C$, and fix once and for all a local coordinate $z$
at $p$.  More precisely, let $\{U_i \}$ be an open cover of $C$,
and assume that $p \in U_0$, that $p \notin U_i$ for $i \ne 0$,
and that $z \in \OO_C(U_0)$ is a coordinate centered at $p$.  A
calculation then shows: 

\begin{lem}\label{canon} 
For $a \in \CX$, let $\xi \in H^1(\OO_C)$
be the image of $a/z$ under the coboundary map induced by the
short exact sequence \[0 \lra \OO_C \lra \OO_C(p) \lra \OO_C(p)|_p
\lra 0.\] Let $s\in H^0(L)$ be a section such that $s(p)=0$.  Then
$s \cup \xi =0$ in $H^1(L)$, and in fact, choosing $\xi$ to be
given by the $1$-cocycle \[ \xi_{ij}= \left\{ \begin{array}{ll}
a/z, &\textnormal{if } i=0;\\ 0, &\textnormal{if } i\ne 0,
\end{array}\right. \] and $\ss^{(1)}$ to be the $0$-cochain defined by
\[ \ss^{(1)}_{i}= \left\{ \begin{array}{ll} as/z, &\textnormal{if }
i=0;\\ 0, &\textnormal{if } i\ne0, \end{array}\right. \] then $s
\cup \xi = \delta\ss^{(1)}$, where $\delta$ is the $\check{C}$ech
coboundary map. In other words, if $\LL_1$ is the first order
deformation of $L$ with transition functions
$\ll_{ij}(1+\xi_{ij}t)$, then $s_i+\ss^{(1)}_{i}t$ is a lifting of $s$
to first order. \QED 
\end{lem}

\subsection{The line bundles $\LL_{p;a}^{\pm}$}

For the rest of the paper, we will focus on a particular class of
deformations. For a fixed point $p \in C$, consider the line
bundles $\pi_1^*\OO_{C}(p) \otimes \OO_{{C}\times {C}}(-\Delta)$,
and $\pi_1^*\OO_{C}(-p) \otimes \OO_{{C}\times {C}}(\Delta)$, over
$C\times C$, where $\Delta \subseteq {C} \times {C}$ is the
diagonal.  We fix the coordinate $z$ centered at $p$ as before,
and let $t$ be the coordinate $z$, viewed as a coordinate on an
affine open subset of the second copy of ${C}$. Let $S\subseteq
\CX$ be a small disk, and $u:S\ra C$ be the inverse to $t$.  On $C
\times S$, define the line bundles $\Lambda_{p}^{\pm}=(1\times
u)^*( \pi_1^*\OO_{C}(\pm p) \otimes \OO_{{C}\times {C}}(\mp
\Delta))$.  From the definition, it is clear that
$\Lambda_p^-=(\Lambda_p^+)^{-1}$. These line bundles induce
holomorphic maps $f:S \ra JC$, and if $w \in S$, then
$(\Lambda^{\pm}_p)_{w}=\OO_{C}(\pm p \mp u(w))$. Using the open
cover of $C$ as in Lemma \ref{canon}, it follows that the
transition functions for $\Lambda^+_p$ are equal to $1$, if
neither $i$ nor $j$ is zero, and (for small $t$), \[ \ll_{0j}(t)=
\frac{z}{z-t} = \sum_{k=0}^{\infty}\left( \frac{t^k}{z^k}
\right).\] Similarly, the transition functions for $\Lambda^-_p$
are equal to $1$, if neither $i$ nor $j$ is zero, and (for small
$t$), \[ \ll_{0j}(t)= \frac{z-t}{z}= 1- \frac{t}{z}.\]

For $L \in \PIC^d(\CT)$, define $\LL^\pm_p= \Lambda_p^\pm \otimes
\pi_1^*L$, a line bundle over $C\times S$. If the transition
functions for $L$ are given by $\ll_{ij}$, then it follows that
the transition functions for $\LL^+_p$ are equal to $\ll_{ij}$, if
neither $i$ nor $j$ is zero, and (for small $t$), \[ \ll_{0j}(t)=
\ll_{0j} \cdot \left( \frac{z}{z-t} \right)= \ll_{0j}\cdot
\sum_{k=0}^{\infty}\left( \frac{t^k}{z^k} \right).\] Similarly,
the transition functions for $\LL^-_p$ are equal to $\ll_{ij}$, if
neither $i$ nor $j$ is zero, and (for small $t$), \[ \ll_{0j}(t)=
\ll_{0j} \cdot \left( \frac{z-t}{z} \right)= \ll_{0j}\cdot \left(
1- \frac{t}{z}\right).\]

For making computations, it will be useful to rescale $t$. For
$a\in \CX$, define a local deformation $\LL^+_{p;a}$ by setting
the transition functions equal to $\ll_{ij}$, if neither $i$ nor
$j$ is zero, and (for small $t$), \[ \ll_{0j}(t)= \ll_{0j} \cdot
\left( \frac{z}{z-at} \right)= \ll_{0j}\cdot
\sum_{k=0}^{\infty}\left( \frac{a^k}{z^k}t^k \right).\] In short
we are considering a second small disk $S'$ in $\CX$, a map $S'\ra
S$ given by $w\mapsto aw$, and setting
$\LL_{p;a}^+=\LL^+_p|_{C\times S'}$.  There is then an induced
holomorphic map $f_a:S' \ra JC$ for each $a$, and if $w \in S'$,
then $(\LL^+_{p;a})_{w}=L\otimes \OO_{C}( p - u(aw))$.

Similarly, define a local deformation $\LL^-_{p;a}$ by setting the
transition functions equal to $\ll_{ij}$, if neither $i$ nor $j$
is zero, and (for small $t$), \[ \ll_{0j}(t)= \ll_{0j} \cdot
\left( \frac{z-at}{z} \right)= \ll_{0j}\cdot \left( 1-
\frac{a}{z}t\right).\] If $w \in S'$, then
$(\LL^-_{p;a})_{w}=L\otimes \OO_{C}( -p + u(aw))$.

A section $s \in H^0(L)$ lifts to first order as a section of
$(\LL^+_{p;a})_1$ if and only if there exists $\ss^{(1)}$
satisfying \[\ss_i^{(1)}-\ll_{ij}\ss_{j}^{(1)}=
\ll_{ij}\aa_{ij}^{(1)}s_j= \left\{ \begin{array}{ll}
 0,& \textnormal{for }i \ne 0;\\
a s_0/z,& \textnormal{for }i=0,
\end{array} \right. \]
and a section $(s+\ss^{(1)}t) \in H^0((\LL^+_{p;a})_1)$
lifts to second order if and only if there exists
$\ss^{(2)}$ satisfying
{\small\[\ss_i^{(2)}-\ll_{ij}\ss_{j}^{(2)}=
\ll_{ij}\aa_{ij}^{(2)}s_j+\ll_{ij}\aa_{ij}^{(1)}\ss_j^{(1)}=\]
\[
\left\{ \begin{array}{ll}
 0,\ \ \ \  \textnormal{for }i \ne 0;& \\
a^2 s_0/ z^2+
\ll_{ij}a \ss^{(1)}_j/z=
a^2 s_0/z^2+
\left(a \ss^{(1)}_0 /z-a^2 s_0/z^2\right)=
a \ss_0^{(1)}/z,& \textnormal{for }i=0.
\end{array} \right.\]}

Likewise, a section $s \in H^0(L)$ lifts to first order as a section 
of $(\LL^-_{p;a})_1$ if and only if there exists $\ss^{(1)}$
satisfying
\[\ss_i^{(1)}-\ll_{ij}\ss_{j}^{(1)}=
\ll_{ij}\aa_{ij}^{(1)}s_j=
\left\{ \begin{array}{ll}
 0,& \textnormal{for }i \ne 0;\\
-a s_0/z,&\textnormal{for }i=0,
\end{array} \right. \]
and a section $(s+\ss^{(1)}t) \in H^0((\LL^-_{p;a})_1)$
lifts to second order if and only if there exists
$\ss^{(2)}$ satisfying
\[\ss_i^{(2)}-\ll_{ij}\ss_{j}^{(2)}=\]
\[
\ll_{ij}\aa_{ij}^{(2)}s_j+\ll_{ij}\aa_{ij}^{(1)}\ss_j^{(1)}=
\left\{ \begin{array}{ll}
 0,\ \ \ \  \textnormal{for }i \ne 0;& \\
\ll_{ij}\left( - a/z\right)\ss^{(1)}_j=
-a \ss^{(1)}_0/z-a^2 s_0/z^2,&
\textnormal{for }i=0.
\end{array} \right.\]

A straightforward calculation 
in the $\CH$ech complex will prove the following two lemmas:

\begin{lem}\label{n1n21}
    Let $\partial_{L,p}$ be the coboundary map 
    $$
    \partial_{L,p}: H^{0}(L(p)\otimes \OO_{p}) \lra H^{1}(L)
    $$
    induced from the exact sequence
    $$
    0\lra L \lra L(p) \lra L(p) \otimes \OO_{p} \lra 0,
    $$
    and let $A_{1}(s)\in H^{0}(L(p) \otimes \OO_{p})$ be defined as
    \[ A_1(s) =
\left\{ \begin{array}{ll}
-as_0/z, & \textnormal{ for } \LL^+_{p;a};\\
as_0/z,& \textnormal{ for } \LL^-_{p;a}.
\end{array} \right.\]
    Then $\gamma_{1}(s)=\partial_{L,p}(A_{1}(s))$.\QED
\end{lem}

\begin{cor}
If $s \in H^0(L(-p))$, then $\gamma_1(s)=0$. \QED
\end{cor}
\remskip

\begin{rem} This is a weaker statement than was proven in Lemma
\ref{canon}, where an explicit first order lifting of $s$ was
given. \end{rem}

\begin{lem}\label{n1n22}
    Let $\partial_{L,2p}$ be the coboundary map 
    $$
    \partial_{L,2p}: H^{0}(L(2p)\otimes \OO_{2p}) \lra H^{1}(L)
    $$
    induced from the exact sequence
    $$
    0\lra L \lra L(2p) \lra L(2p) \otimes \OO_{2p} \lra 0,
    $$
    and let $A_{2}(s+\sigma^{(1)}t)
\in H^{0}(L(2p) \otimes \OO_{2p})$ be defined as
    \[ A_2(s+\ss^{(1)}t) =
\left\{ \begin{array}{ll}
-a\ss^{(1)}_0/z, & \textnormal{ for } \LL^+_{p;a};\\
a\ss^{(1)}_0/z+a^2s_0/z^2,& \textnormal{ for } \LL^-_{p;a}.
\end{array} \right.\]
    Then 
    $\gamma_{2}(s+\sigma^{(1)}t)=\partial_{L,2p}(A_{2}(s+\sigma^{(1)}t))$.\QED
\end{lem}

\begin{cor}\label{seccor} If $s \in H^0(L(-2p))$, then there
exists a first order lift of $s$, say $s +\ss^{(1)}t$, such that
$\gamma_2(s+\ss^{(1)}t)=0$. \end{cor}

\PRF.  Let $s \in H^{0}(L(-2p))$.  Then since $s \in
H^{0}(L(-p))$, let $s+\ss^{(1)}t$ be the standard lift of $s$, as
given in Lemma \ref{canon}.  Recall that we set
\begin{displaymath} \sigma^{(1)}_i= \left\{ \begin{array}{ll}
    as/z, & \textnormal{if } i=0,\\
   0, & \textnormal{if } i \ne 0.
   \end{array} \right.
 \end{displaymath} 
Since $s$ vanishes to order 2 at $p$, we see that
$A_{2}(s+\sigma^{(1)})=0$, and hence $s$ lifts to second order.
\QED \remskip

For our computations, we will want to consider a more general
class of deformations modeled on the $\LL_{p;a}^\pm$.  Let
$\Delta_i\subseteq C\times C^k$ be defined as
$\Delta_i=\{(x_0,\ldots,x_k) \in C\times C^k |x_0=x_i\}$, and for
a particular choice of points $p_1,\ldots,p_k$, let
$D=\sum_{i=1}^kp_i$.  On $C \times C^k$, consider the line bundle
$\pi_1^*\OO_C(D)\otimes \OO_{C\times C^k}(-\sum\Delta_i)$. Letting
$u_i$ be a map from a disk $S\subseteq \CX$ to a neighborhood of
the point $p_i$, and $u:S\ra C^k$ be the map given by $w \mapsto
(u_1(w),\ldots,u_k(w))$, then set $\Lambda^+_{D}=(1\times
u)^*(\pi_1^*\OO_C(D)\otimes \OO_{C\times C^k}(-\sum\Delta_i))$. 
For $L\in \PIC ^d(C)$, let $\LL_D^+=\Lambda_D^+ \otimes \pi_1^*L$.
This has fiber over a point $w\in S$ equal to $L\otimes
\OO_{C}((p_1-u_1(w)) \otimes \ldots \otimes \OO_{C}((p_k-u_k(w))$. 
As before $\Lambda_D^+$ induces a holomorphic map $S\ra JC$, and
in addition it is clear that $\LL^+_D= \pi_1^*L \otimes
\Lambda^+_{p_1}\otimes \ldots \otimes \Lambda^+_{p_k}$.  We can
similarly define $\LL^{-}_D$, and by rescaling the local
coordinate, construct the line bundle
$$
\LL= \pi_1^*L \otimes 
\Lambda^+_{p_1;a_1}\otimes \ldots \otimes \Lambda^+_{p_{k_1};a_{k_1}} \otimes
\Lambda^-_{p_{k_{1}+1};a_{k_{1}+1}}\otimes 
\ldots \otimes \Lambda^-_{p_{k};a_{k}}.
$$
The fiber of $\LL$ over a point $w\in S$ is given by
$$
\LL_{w} = L \otimes \OO_{C}(p_1-u_1(a_1w))
\otimes \ldots 
\otimes \OO_{C}(-p_{k}+u_k(a_{k}w)).
$$

\begin{rem} The calculations in Lemmas \ref{n1n21} and \ref{n1n22}
are local, in the sense that the obstructions for $\LL$ are sums
of the local contributions calculated in those lemmas.  \end{rem}

\subsection{The Prym case}

For the rest of the paper, we will be considering the following
situation.  $C$ will be a smooth curve of genus $g$, $\pi:\CT \ra
C$ will be a connected \'etale double cover, $\tt$ will be the
associated involution on $\CT$, $\eta\in \PIC^0(C)$ will be the
associated semiperiod, and $P \subseteq J\CT$ will be the Prym
variety. If $\widetilde{\Theta}$ is the canonical theta divisor of
$J\CT$, then Mumford \cite{mum1} has shown that
$\widetilde{\Theta} \cap P= 2 \cdot \Xi$, where $\Xi$ is the class
of a principal polarization on $P$.  Recall, if we identify $JC$
with $\PIC^{2g-2}(\CT)$, then as a set $P$ can be described as
$$P= \{ L \in \PIC^{2g-2}(\CT) \ | \
\textnormal{Norm}(L)=\omega_C, \ h^0(L) \equiv 0 \
(\textnormal{mod } 2) \}.$$ It follows that $\Xi=\{ L \in P \ | \
h^0(L) \ge 2 \}$. 

The following straightforward lemma is the fundamental tool that
we will use in what follows. 

\begin{lem}\label{foundlem}
Let $H$ be a hypersurface, not necessarily reduced, defined in an
open neighborhood of $0$ in $\CX^{n}$ and containing $0$.  Let $S$
be a disk in $\CX$, containing $0$, and $f:S \ra \CX^{n}$ be a
holomorphic map, with $f(0)=0$.  Then $\MUL_{0}H \le
\MUL_{0}f^{*}H$, and equality holds if and only if $f_{*}(T_{0}S)$
is not contained in the tangent cone to $H$ at $0$. 
\hspace*{\stretch{1}}$\square$ 
\end{lem}

As an application, suppose $S$ is a smooth curve with $s_0 \in S$
and $\LL$ is a line bundle over $\CT \times S$ of relative degree
$2g-2$.  Let $f:S \ra J\CT$ be the induced morphism, 
and let $L\in \PIC^{2g-2}(\CT)$ be the line bundle associated to the point
$x=f(s_0) \in J\CT$.

\begin{lem}[\cite{cmf}, 1.4]\label{lemprymmult}  
If $f(S)\subseteq P$, then \[ \frac{1}{2}h^0(L) \le \MUL_x\Xi \le
\frac{1}{2}\deg_{s_0} \Theta_S =
\frac{1}{2}\ell((R^1\pi_{2*}\LL)_{s_0}).\] Moreover, there exists
a choice of $S$ and a line bundle $\LL$ as above such that
$\MUL_x\Xi=\frac{1}{2}\ell((R^1\pi_{2*}\LL)_{s_0})$. \QED
\end{lem}

Recall that for $t$ a local coordinate on $S$ centered at $s_{0}$
and only vanishing there, $\CT_{k}= \CT \times \SPEC
\CX[t]/(t^{k+1})$, $W_k$ is the image of the map $H^0(\LL_k) \ra
H^0(L)$ induced by the exact sequence (\ref{sesk}), and
$d_k=\dim(W_k)$. 

\begin{pro}\label{length}
In the notation above, 
\begin{itemize}
    \item[\textnormal{(a)}]  For every $k$,
$\ell((R^1\pi_{2*}\LL)_{s_0})\ge \ell(H^{0}(\LL_{k}))= 
        \sum_{i=0}^{k}d_{i}$, and if $d_N=0$, then 
equality holds for all $k\ge N$;
        
\item[\textnormal{(b)}]  if $f(S)\subseteq P$, and 
$L= \pi^*(M)\otimes \OO_{\CT}(B)$, with $h^0(M)> h^0(L)/2$, $B\ge 0$,
and $B\cap \tt^*B=\emptyset$, then $d_1 \ge 2h^0(M)-h^0(L)$;

\item[\textnormal{(c)}] 
if $x \in \Xi$, and $L= \pi^*(M)\otimes \OO_{\CT}(B)$, with $h^0(M)>
h^0(L)/2$,  $B\ge 0$, and $B\cap \tt^*B=\emptyset$, then 
$\MUL_x \Xi \ge  h^0(M)$.  Furthermore, $\MUL_x \Xi = h^0(M)$ if and only if
there exists a choice of $S$ and a line bundle $\LL$ as above such that
$f(S)\subseteq P$, $d_1=2h^0(M)-h^0(L)$, and $d_2=0$.

\end{itemize}   
\end{pro}

\PRF.  (a) is a restatement of the lemmas in Section 1.1, and
the proof of (b) is contained in the proof of 
\cite{cmf} Theorem 2.3.

(c) It follows from part (a) and Lemma \ref{lemprymmult} that for 
any deformation  $f:S \ra P$ with $f(s_0)=x$,

$$
\MUL_x \Xi \le \frac{1}{2} \MUL_{s_0}\Theta_S =
\frac{1}{2}\sum_{k\ge 0} d_k. $$ Furthermore, there exist
deformations for which equality holds, so that $$\MUL_x \Xi =
\inf_{f:S \ra P; f(s_0)=x}\{\frac{1}{2}\sum_{k\ge 0} d_k\}. $$ By
definition  $d_0=h^0(L)$, and by (b) we know that $d_1\ge 2h^0(M)
-h^0(L)$ for all such deformations. Since $d_k \ge d_{k+1}\ge 0$
for all $k$, it follows that $$\MUL_x \Xi = \inf
\{\frac{1}{2}\sum_{k\ge 0} d_k\} \ge \frac{1}{2}
(h^0(L)+2h^0(M)-h^0(L))= h^0(M).$$ Equality holds if and only if
there is a deformation such that $d_1=2h^0(M)-h^0(L)$, and
$d_2=0$.  \QED

\remskip

Let $\{U_i \}_{i\in I}$ be an open affine cover for $\CT$, where
for $i\in \{1,\ldots,n\}$, $\ p_{i} \in U_{j} \ \iff \ i=j$, and
$\tau(p_{i}) \in U_{j} \ \iff \ j=i+n$.  For $ 1\le i\le n$ we will
define the index $\tau(i)=i +n$.  On each open set $U_{i}$ define
$z_{i}$ to be a local coordinate, which for $i\in \{1,\ldots,n\}$
is centered at $p_{i}$, and for $i\in \{\tau(1),\ldots,\tau(n)\}$
is centered at $\tau(p_{i})$. In addition, we will choose the
local coordinates so that $\tau^*z_i=z_{\tau(i)}$. 
  
Let $q_{1},\ldots,q_n$ be general points of $C$, $\pi^{-1}(q_{i})=
\{p_{i},\tau(p_{i})\}$, and $ D =
\sum_{i=1}^{n}(p_{i}+\tau(p_{i}))$.  Let
$a=(a_{1},\ldots,a_{n})\in \CX^{n}$. \remskip

\begin{dfn}\label{gooddef} Let
$S\subseteq \CX$ be a disk containing the origin, and let $L\in
\PIC^{2g-2}(\CT)$ be the line bundle associated to a point $x\in
P$.  With $D$ and $a$ as above, 
define the \emph{deformation of $L$ associated to $D$ and
$a$}, denoted by $\LL_{D;a}=\LL$, to be the line bundle over $\CT\times S$ 
given by $$\LL=\pi_1^*L\otimes \Lambda^+_{p_1;a_1}\otimes \ldots
\otimes \Lambda^+_{p_n;a_{n}} \otimes
\Lambda^-_{\tt(p_1);a_{1}}\otimes
 \ldots \otimes \Lambda^-_{\tt(p_{n});a_{n}}.$$
The fiber over a point $w\in S$ is given by
\[ \LL_w=L \otimes \OO_{\CT}(p_1-u_1(a_1w)-\tau(p_1)+\tau(u_1(a_1w)))
\otimes \ldots \]
\[\otimes \OO_{\CT}\left(p_n-u_n(a_nw)-\tau(p_n)+\tau(u_n(a_nw))\right).\]
\end{dfn}

\begin{lem}
Given $\LL_{D,a}$, let $f:S\ra J\CT$ be the associated morphism.
Then $f$ is holomorphic, and $f(S)\subseteq P$. 
\end{lem}

\PRF. We have seen in the previous section that $f$ is
holomorphic. Let $s_0\in S$ be such that $f(s_0)=x$.  
Then $f(S) \subseteq P$, since $\textnormal{Norm}((\LL_{D;a})_s)= \omega_C$ 
for all $s \in S$, and
$(\LL_{D;a})_{s_0}=L\in P$. Indeed, $\textnormal{Norm}^{-1}(\omega_C)$
has two connected components distinguished by the parity of
$h^0$. Since $f(S)$ includes a point in the Prym variety, namely
$L$, $f(S)$ is contained in $P$.\QED \remskip

We will now reinterpret Lemmas \ref{n1n21} and \ref{n1n22} in the
case of a deformation $\LL_{D;a}$.  To begin, we fix a
trivialization of $L$ and $\LL_1$ at the points
$p_1,\tt(p_1),\ldots, p_n,\tt(p_n)$.  We then choose once and for
all a basis for $H^0(L(D) \otimes \OO_D) $ and $H^0(L(2D) \otimes
\OO_{2D}) $; in the former case we will take
$\{1/z_1,1/z_{\tt(1)},\ldots,$ $1/z_n,1/z_{\tt(n)} \} $, and in
the latter case $\{1/z_1^2,1/z_1,1/z_{\tt(1)}^2,$
$1/z_{\tt(1)},\ldots, 1/z_n^2,$ $1/z_n,1/z_{\tt(n)}^2,1/z_{\tt(n)}
\} $.  With respect to these bases and trivializations, the
lemmas can then be restated as follows: 

\begin{lem}\label{Dn1n21}
 Let $\partial_{L,D}$ be the coboundary map 
    $$
    \partial_{L,D}: H^{0}(L(D)\otimes \OO_{D}) \lra H^{1}(L)
    $$
    induced from the exact sequence
    $$
    0\lra L \lra L(D) \lra L(D) \otimes \OO_{D} \lra 0,
    $$
    and let $A_{1}(s)\in H^{0}(L(D) \otimes \OO_{D})$ be defined as
    \[ A_1(s) =
(-a_1s(p_1),a_1s(\tt(p_1)),\ldots,
-a_ns(p_n), a_ns(\tt(p_n)).\]

Then $s\in H^0(L)$ lifts to first order as a section of 
$\LL_{D;a}$ if and only if 
$s\in \ker(\partial_{L,D}\circ A_1)$.\QED

\end{lem}

\begin{cor}
$H^0(L(-D))\subseteq \ker(A_1)$, and if $a_i\ne 0$ for all $i$, then 
$H^0(L(-D))=\ker(A_1)$.\QED
\end{cor}

\begin{lem}\label{Dn1n22}
 Let $\partial_{L,2D}$ be the coboundary map 
    $$
    \partial_{L,2D}: H^{0}(L(2D)\otimes \OO_{2D}) \lra H^{1}(L)
    $$
    induced from the exact sequence
    $$
    0\lra L \lra L(2D) \lra L(2D) \otimes \OO_{2D} \lra 0,
    $$
    and let $A_{2}(s+\ss^{(1)}t)\in H^{0}(L(2D) \otimes \OO_{2D})$ 
be defined as
{\small
$$
A_{2}(s+\sigma^{(1)}t)=(0,-a_{1}\sigma^{(1)}(p_{1}),\ a_{1}^{2}s(\tau (p_{1})),
\ a_{1}\sigma^{(1)}(\tau (p_{1})) +a_{1}^{2}\frac{ds}{dz}(\tau 
(p_{1})),\ldots).
$$}
Then $s+\ss^{(1)}t \in H^0(\LL_1)$ lifts to second 
order 
if and only if 
$s\in \ker(\partial_{L,2D}\circ A_2)$.\QED

\end{lem}

\begin{cor}\label{cor149}
$H^0(L(-2D))\subseteq W_2$.
\end{cor}

\PRF.  $s\in W_2$ if and only if there is a first
order lift $s+\ss^{(1)}t$ such that $\partial_{L,2D}\circ A_2
(s+\ss^{(1)}t)=0$. Using the standard lift of $s$ given in Lemma
\ref{canon}, and the same analysis as in Corollary \ref{seccor},
one can easily show that if $s\in H^0(L(-2D))$, then $A_2
(s+\ss^{(1)}t)=0$.\QED

\section{Linear Systems and Linear Algebra}

In this section we collect some general results that will be
useful for making computations in subsequent sections.

\subsection{Linear systems on a double cover}

Given a line bundle $L$ on $\CT$, and a point $q\in C$, we will
want to know when the points of $\pi^{-1}(q)$ impose independent
conditions on the linear system $|L|$. 

\begin{lem}\label{conditions} Suppose that $L=\pi^*M \otimes
\OO_{\CT}(B)$, where $M$ is a line bundle on $C$ such that
$h^0(\CT,L)>h^0(C,M)>0$, and $B \ge0$ is an effective divisor on
$\CT$ such that $B \cap \tau^*B= \emptyset$.  If $p\in \CT$ is a
general point, then $h^0(L(-p-\tt(p)))=h^0(L)-2$. 

\end{lem}

\PRF.  Let $b\in H^0(\OO_{\CT}(B))$ be such that $(b)_0=B$, and
let $v\in H^0(L) - \pi^*H^0(M)\cdot b$.  We can consider
$(\frac{v}{b})$ as a rational section of $\pi^*M$, and I claim
that $(\frac{v}{b})$ is not $\tt$-invariant.  Indeed, if it were,
then it could not have poles along $B$ since it has none along
$\tt(B)$, and so $(\frac{v}{b}) $ would be a regular section; i.e. 
$\frac{v}{b} \in H^0(\pi^*M)^+=\pi^*H^0(M)$.  This would be a
contradiction, as then $v \in \pi^*H^0(M)\cdot b$.

It follows that for a general point $p\in \CT$,
$(\frac{v}{b})(p)\ne (\frac{v}{b})(\tt(p))$.  Let $s \in
\pi^*H^0(M)$ be a nonzero section, and let $\ll\in \CX$ be such
that $v(p) - \ll s(p) b(p) =0$.  It is immediate to check that
$v(\tt(p))-\ll s(\tt(p))b(\tt(p))\ne 0$.  This completes the
proof, as for a general $p$, $h^0(L(-p))=h^0(L)-1$, and we have
exhibited a section $v-\ll s b \in H^0(L(-p))$ which does not
vanish at $\tt(p)$.\QED

\remskip 
\begin{cor}\label{cor212} Suppose that $L=\pi^*M \otimes
\OO_{\CT}(B)$, where $M$ is a line bundle on $C$ such that
$h^0(\CT,L)\ge h^0(C,M)>0$, and $B \ge0$ is an effective divisor
on $\CT$ such that $B \cap \tau^*B= \emptyset$.  Let $b\in H^0(B)$
be a section vanishing on $B$. Then $H^0(L)=\pi^*H^0(M)\cdot b$ if
and only if for a general point $p\in \CT$
$h^0(L(-p-\tt(p)))=h^0(L)-1$. \end{cor}

\PRF. Let $p$ be a general point of $\CT$, and $q=\pi(p)$.  If
$H^0(L)=\pi^*H^0(M)\cdot b$, then it follows that
$H^0(L(-p-\tt(p)))=\pi^*H^0(M(-q))\cdot b$.  Hence,
$h^0(L(-p-\tt(p)))=h^0(L)-1$. Conversely, if for a general point
$p$, $h^0(L(-p-\tt(p)))=h^0(L)-1$, then the previous lemma implies
that $h^0(L)=h^0(M)$. \QED

\remskip 
For the duration of the paper we will use the following
notation. Given a collection $p_1,\tt(p_1),\ldots,p_k,\tt(p_k)$ of distinct
points of $\CT$, we will set $D_k=\sum_{i=1}^{k}(p_i+\tt(p_i))$.

\begin{cor}\label{cor313}

Suppose that $L=\pi^*M \otimes \OO_{\CT}(B)$, where $M$ is a line
bundle on $C$ such that $h^0(\CT,L)\ge h^0(C,M)>0$, and $B \ge0$
is an effective divisor on $\CT$ such that $B \cap \tau^*B=
\emptyset$.  Let $h^0(M)=n_1$, $h^0(L)-h^0(M)=n_2$, and
$p_1,\tt(p_1),\ldots,p_{k},\tt(p_{k})$ be $2k$ points of $\CT$,
where $p_1,\ldots,p_k$ are general.

\begin{itemize}
\item[\textnormal{(a)}] Suppose $h^0(M)>h^0(L)/2$ and 
$k\le n_2$. 
Then
$h^0(L(-D_k))=h^0(L)-2k$.

\item[\textnormal{(b)}] Suppose $h^0(M)>h^0(L)/2$ and 
$ n_2 \le k \le n_1$.  Then 
$h^0(L(-D_k))=h^0(L)-n_2-k$.
Furthermore, in this case
$$H^0(L(-D_k))=
\pi^*H^0(M(-\sum_{i=1}^{k}q_i))\cdot b,$$
where $q_i=\pi(p_i)$ and $b\in H^0(\OO_{\CT}(B))$ vanishes on $B$.

\item[\textnormal{(c)}]
 Suppose $h^0(M)>h^0(L)/2$, $ n_2 \le k \le n_1$, and $1\le k_1 \le k$.
Then 
{\small $H^0(L(-D_k-D_{k_1}))=
H^0(L(-D_k-\sum_{i=1}^{k_1}p_i))=
H^0(L(-D_k-\sum_{i=1}^{k_1}\tt(p_i)))$}
and 
$h^0(L(-D_k-D_{k_1}))=\max (h^0(L)-n_2-k-k_1,0)$.
Furthermore, in this case
$$H^0(L(-D_k-D_{k_1}))=
\pi^*H^0(M(-\sum_{i=1}^{k}q_i-\sum_{j=1}^{k_1}q_j))\cdot b.$$

\item[\textnormal{(d)}] Suppose $h^0(M)\le h^0(L)/2$
and $k\le n_1$.
Then
$h^0(L(-D_k))
=h^0(L)-2k$.
\end{itemize}
\end{cor}

\PRF.  (a) In the case $h^0(L)=h^0(M)$, there is nothing to prove
since $k \le n_2=0$.  So assume that $h^0(L)> h^0(M)$.  We will
now use induction on $k$.  For $k=1$, we are done by the previous
lemma. So assume we have proven (a) for all $k \le m $, where $m$
is some integer less than $n_2$.  Let $q_i=\pi(p_i)$,
$L'=L(-D_m)$, and $M'=M(-\sum_{i=1}^{m}q_i)$.  Then
$L'=\pi^*M'\otimes \OO_{\CT}(B)$, $h^0(L')=h^0(L)-2m$ by
induction, and $h^0(M')=h^0(M)-m$ since the points $q_i$ are
general. Now $h^0(L')=h^0(L)-2m > h^0(M)-m=h^0(M')>0$, since
$m<n_2<h^0(M)$ and $h^0(L)-h^0(M)=n_2 >m$.  It follows that $L'$
and $M'$ satisfy the condition of the lemma, and hence for general
points $ p_{m+1}$ and $\tt(p_{m+1})$ we have
$h^0(L(-D_{m+1}))=h^0(L)-2m-2$. 

(b) In the case $k=n_2$, we have seen that
$h^0(L(-D_{n_2}))=h^0(L)-2n_2=n_1-n_2
=h^0(M(-\sum_{i=1}^{n_2}q_i))$.  Hence the natural inclusion $$
\pi^*H^0(M(-\sum_{i=1}^{n_2}q_i)))\cdot b\subseteq
H^0(L(-D_{n_2})) $$ is an equality.  For $k>n_2$ we use induction
and the previous corollary. 

(c) By part (b), $H^0(L(-D_k))=
\pi^*H^0(M(-\sum_{i=1}^{k}q_i))\cdot b$.  Thus the vanishing locus
of a section of $H^0(L(-D_k))$ is invariant away from the support
of $B$.  It follows that
$H^0(L(-D_k-p_1-\tt(p_1)))=H^0(L(-D_k-p_1))=H^0(L(-D_k-\tt(p_1)))
=\pi^*H^0(M(-\sum_{i=1}^{k}q_i-q_1))\cdot b$, since the points
$q_i$ are general.  We also have that
$h^0(M(-\sum_{i=1}^{k}q_i-q_1))=
\max(h^0(M(-\sum_{i=1}^{k}q_i))-1,0)$, since if $q_1$ is a base
point of $H^0(M(-\sum_{i=1}^{k}q_i))$, then it is a ramification
point of the map associated to $|M|$.  (The ramification locus is
finite, and the $q_i$ are general, so we can assume $q_1$ is not a
base point of $H^0(M(-\sum_{i=1}^{k}q_i))$.) One can then proceed
by induction on $k_1$. 

(d) An induction argument similar to that in part (a) will prove
this.\QED

\begin{lem}\label{shpair}
Suppose $h^0(L)=2n>0$, and that for every choice of distinct
points 
$p_1,\tt(p_1),\ldots,p_n,\tt(p_n)$ of $\CT$, $h^0(L(-D_n))>0$.
Then $L=\pi^*M \otimes \OO_{\CT}(B)$, where $M$ is a line bundle
on $C$ such that $h^0(C,M)>h^0(\CT,L)/2$, $B \ge0$ is an effective
divisor on $\CT$ such that $B \cap \tau^*B= \emptyset$, and
$h^0(\CT,B)=1$.  \end{lem}

\PRF.  For every $s\in H^0(L)$, the divisor $(s)_0$ can be
decomposed into an invariant part say $N=\pi^*N'$, and the
residual part, say $B$, which by definition must have the property
$B\cap\tau^*B=\emptyset$. Hence, setting $M=\OO_C(N')$ we can
always write $L=\pi^*M\otimes \OO_{\CT}(B)$, with $h^0(L)\ge
h^0(M) >0$, $B \cap \tau^*B=\emptyset$, and $B\ge 0$. 

We will first prove the lemma in the case that there do not exist
$2$ points $p$ and $\tau(p)$ such that
$h^0(L(-p-\tt(p)))=h^0(L)-2$. In this case, Corollary \ref{cor212}
implies that $h^0(M)=h^0(L)$.  Since $B$ is effective, $h^0(B)\ge
1$, and the inequality $\dim|M| = \dim|L|\ge \dim |\pi^*M|
+\dim|B| \ge \dim|M|+\dim|B|$ implies $h^0(B)=1$. Hence $L=\pi^*M
\otimes \OO_{\CT}(B)$, $h^0(L)=h^0(M)$, $B \cap \tt^*B =
\emptyset$, and $h^0(B)=1$. 

We will now proceed to prove the lemma by induction on $h^0(L)$.
The case $h^0(L)=2$ is a consequence of the case above. So suppose
we have proven the result for all line bundles $L'$ for which
$h^0(L')\le 2n-2$, and consider a line bundle $L$ with
$h^0(L)=2n>2$.  By the case above, we may assume there exist
points $p$ and $\tt(p)$ imposing independent conditions on
$H^0(L)$.  Let $L'=L(-p-\tt(p))$, so that $h^0(L')=2n-2$.  There
do not exist $2n-2$ distinct points $p_1,\tt(p_1),\ldots
p_{n-1},\tt(p_{n-1})$ on $\CT$ imposing independent conditions on
$H^0(L')$, since otherwise, after possibly replacing
$p_1,\tt(p_1),\ldots
p_{n-1},\tt(p_{n-1})$ with
 a more general choice of points,  $p,\tt(p),p_1,\tt(p_1),\ldots
p_{n-1},\tt(p_{n-1})$ would be distinct points imposing independent conditions on
$H^0(L)$, contradicting our assumptions.

Thus, by induction, $L'=\pi^*M'\otimes \OO_{\CT}(B)$, with
$h^0(M')>n-1$, $B\cap \tt^*B =\emptyset$, $B\ge 0$ and $h^0(B)=1$.
Setting $q=\pi(p)$, it follows that $L=\pi^*(M'(q))\otimes
\OO_{\CT}(B)$.  If we let $M=M'(q)$, then $h^0(M)\ge h^0(M')\ge
n$. In the case $h^0(M)=n$, we would arrive at a contradiction,
since part (c) of the previous corollary with $k=n_1=n$ would
imply that there were $2n$ distinct points
$p_1,\tt(p_1),\ldots,p_n,\tt(p_n)$ imposing independent conditions
on $H^0(L)$. Thus, $h^0(M)>n$. \QED

\begin{lem}[Smith-Varley \cite{sv2}]\label{uniqueM}
Suppose $L=\pi^*M \otimes \OO_{\CT}(B)$, where $M$ is a line
bundle on $C$ such that $h^0(C,M)>h^0(\CT,L)/2>0$, and $B \ge0$ is
an effective divisor on $\CT$ such that $B \cap \tau^*B=
\emptyset$.  Then $M$ and $\OO_{\CT}(B)$ are unique up to
isomorphism. \end{lem}

\PRF.  Suppose $L=\pi^*M' \otimes \OO_{\CT}(B')$, where $M'$ is a
line bundle on $C$ such that $h^0(C,M')>h^0(\CT,L)/2>0$, and $B'
\ge0$ is an effective divisor on $\CT$ such that $B' \cap
\tau^*B'= \emptyset$.  Then $\pi^*H^0(M)\cdot b$ and
$\pi^*H^0(M')\cdot b'$ are both linear subspaces of $H^0(L)$, and
 $\dim(\pi^*H^0(M)\cdot b) + \dim(\pi^*H^0(M')\cdot b')
\ge h^0(L)+2$.  Hence they have nontrivial intersection.
 
Let $s\in \pi^*H^0(M)$ and $s'\in \pi^*H^0(M')$ be sections such
that $s\cdot b=s'\cdot b'$.  Then $(s)_0+(b)_0=(s')_0+(b')_0$. The
invariant parts of these divisors must agree, and so we see
$(s)_0=(s')_0$ and $(b)_0=(b')_0$. Thus $M\cong \OO_C((s)_0)\cong
\OO_C((s')_0)\cong M'$, and $\OO_{\CT}(B)\cong \OO_{\CT}(B')$.\QED

\remskip
\begin{rem}\label{remcomp} Suppose $L$ can be written in the form
$\pi^*M \otimes \OO_{\CT}(B)$, where $M$ is a line bundle on $C$
such that $h^0(C,M)>h^0(\CT,L)/2>0$, and $B \ge0$ is an effective
divisor on $\CT$ such that $B \cap \tau^*B= \emptyset$. Then the
proof of Corollary \ref{cor313} actually shows that if $k$ is the
maximum number such that there exist points
$p_1,\tt(p_1),\ldots,p_k,\tt(p_k)$ imposing $2k$ conditions on
$H^0(L)$, then $h^0(M)=h^0(L)-k$. \end{rem}

\subsection{Subspaces of complementary dimension}

In computations made in subsequent sections we will have to
examine the intersection of linear subspaces of a given vector
space. Specifically, we will be given a vector space $V=\CX^d$,
two fixed subspaces $V_1$ and $V_2$ such that
$\dim(V_1)+\dim(V_2)=d$, and a family $F$ of linear subspaces of
$V$ parameterized by a second copy of $\CX^d$, and defined as
follows: for $a=(a_1,\ldots,a_d)\in \CX^d$, $$ \begin{array}{ccl}
F_a &=& \{v\in V\ | \ \exists\ v^1 \in V_1 \textnormal{ s.t. }
\pi_i(v)=a_i \pi_i(v^1), \ i=1,\ldots,d\},\\ \end{array} $$ where
$\pi_i$ is projection onto the $i$-th factor. Our goal will be to
determine whether or not there exists an $a\in \CX^d$ such that
$F_a \cap V_2 = 0$. 

To begin, let $\dim(V_1)=d_1$ and $\dim(V_2)=d_2$. Define a
coordinate $m$-plane to be the linear subspace of $V$ defined by
the vanishing of $d-m$ of the $\pi_i$. We will now prove the
following proposition.

\begin{pro}\label{linalg}
If the intersection of $V_2$ with each coordinate $d_1$-plane is
trivial, then there is a Zariski open subset $U\subseteq \CX^d$
such that for all $a\in U$, $F_a \cap V_2 = 0$. 
\end{pro}

\PRF.  Let $\{ (v^{1}_{11},\ldots,v^{1}_{1d}),\ldots,
(v^{1}_{d_11},\ldots,v^{1}_{d_1d})\}$ be a basis for $V_1$, and
let $\{ (v^{2}_{11},\ldots,v^{2}_{1d}),\ldots,
(v^{2}_{d_21},\ldots,v^{2}_{d_2d})\}$ be a basis for $V_2$. A
basis for $F_a$ is then given by $\{
(a_1v^{1}_{11},\ldots,a_dv^{1}_{1d}),\ldots,
(a_1v^{1}_{d_11},\ldots,a_dv^{1}_{d_1d})\}$. 

Let $M$ be the following matrix:
\[
M=\left( \begin{array}{rcr}
 a_1v^{1}_{11}  &\ldots&a_dv^{1}_{1d}\\
 \vdots\ \             &     & \vdots\ \ \\
 a_1v^{1}_{d_11}&\ldots&a_dv^{1}_{d_1d}\\
 v^{2}_{11}  &\ldots&v^{2}_{1d}\\
 \vdots \ \            &     & \vdots\ \ \\
 v^{2}_{d_21}&\ldots&v^{2}_{d_2d}\\
\end{array} \right)
\]
It follows that $F_a \cap V_2=0$ if and only if $\det(M)\ne 0$.
We now appeal to the following lemma, where we will use the
notation $\MA(m,n)$ for the space of $m \times n$ matrices over $\CX$.

\begin{lem}\label{det}
    Let $d'<d \in \mathbb{N}$, $A\in \MA(d',d)$ and $B\in \MA(d-d',d)$ 
have
columns 
    $A_{i}$ and $B_{i}$ respectively, and let $C$ be the matrix
    $$C=
\left( \begin{array}{ccc} 
A_1&\ldots &A_d\\
B_1&\ldots &B_d
    \end{array} \right).$$ 
Then,
    $$
\det(C)=    \sum_{i_{1}<\ldots<i_{d'}}(-1)^{\epsilon + i_1+\ldots + i_{d'}}
    \det(A_{i_{1}}\ldots A_{i_{d'}})\det
    (B_{k_{1}}\ldots B_{k_{d-d'}}),
    $$
where $\{ i_{1},\ldots,i_{d'} \} \cup \{ k_{1},\ldots,k_{d-d'} \}
    = \{1,\ldots,d \}$, 
$ k_{1}<\ldots< k_{d-d'} $, 
and $\epsilon$ is an integer satisfying 
 $\epsilon + \frac{d'(d'+1)}{2} \equiv 0
    \ (\textnormal{mod } 2).$
\end{lem}    
\PRF.  Let $D:\MA(d,d)\ra \CX$ be the given by the above formula.
$D(I)=1$, $D$ is alternating, and $D$ is multilinear in the 
columns. Thus $D=\det$. \QED
\remskip

If we let $V^{i}_j=(v^{i}_{1j},\ldots,v^{i}_{d_ij})^T$, 
then as an immediate consequence of the lemma
    $$
\det(M)=    \sum_{i_{1}<\ldots<i_{d_1}}(-1)^{\epsilon + \Sigma i_{j}}
    a_{i_1}\ldots a_{i_{d_1}}
    \det(V^{1}_{i_{1}}\ldots V^{1}_{i_{d_1}})\det
    (V^{2}_{k_{1}}\ldots V^{2}_{k_{d_2}}).
    $$
The monomials in the $a_i$ which appear in the formula above are
distinct. Also, since the dimension of $V_1$ is $d_1$, there must
be some choice of ${i_{1}<\ldots<i_{d_1}}$ such that
$\det(V^{1}_{i_{1}}\ldots V^{1}_{i_{d_1}})\ne 0$. I claim that for
any choice of ${k_{1}<\ldots<k_{d_2}}$, $\det
    (V^{2}_{k_{1}}\ldots V^{2}_{k_{d_2}})\ne 0$.  It follows that
the determinant of $M$ is not identically zero as a polynomial in
the $a_i$, and thus we can take $U = \{ \det(M)\ne 0\}$. 

We now proceed to prove the claim.

\begin{lem}\label{col}
Suppose $\Lambda =\{\ll_{ij}\}\in \MA(d',d)$ has rank $d'$,
and let $H$ be the $d'$ dimensional 
linear subspace of $\CX^d$ spanned by the rows of $\Lambda$.  
Suppose further that one of the following equivalent conditions
hold:
\begin{itemize}

\item[\textnormal{(a)}]
the intersection of $H$ with any coordinate $(d-d')$-plane 
is trivial;

\item[\textnormal{(b)}]
a linear combination of the rows of $\Lambda$ with 
$d'$ or more entries equal to zero is identically zero;

\item[\textnormal{(c)}] 
if 
$\sum_{i=1}^{d'} \aa_i \ll_{ij}=0$ for some $\aa_1,\ldots ,\aa_{d'}
\in \CX$, and
for all $j\in S \subseteq \{1,\ldots,d\}$, with $|S|=d'$,
then $\aa_i=0$ for all $i$.
\end{itemize}
Then every
$d'$ columns of $\Lambda$ are linearly independent.
\end{lem}

\PRF.  If there exist $d'$ columns of $\Lambda$ which are
dependent, then there would be a dependence among the rows of
those columns.  This would imply that some nontrivial linear
combination of the rows of $\Lambda$ had at least $d'$ entries
which were zero, which contradicts our assumption on $\Lambda$.
\QED

\remskip Clearly applying this lemma to the space $V_2$ finishes
the proof of the claim and hence of the proposition. \QED

\remskip Due to Lemma \ref{col}, there is a useful restatement of
the proposition. 

\begin{cor}\label{cor324} In the notation above, suppose $V_2$
satisfies the following condition: if $v_2 \in V_2$, and
$\pi_i(v_2)=0$ for $i\in I\subseteq \{1,\ldots,d\}$ with
$|I|=d_2$, then $v_2=0$.  Then there exists a Zariski open set $U
\in \CX^d$ such that for all $a\in U$, $F_a \cap V_2 = 0$. \QED
\end{cor}

\section{Proof of Theorem \ref{teosv}}

\subsection{Preliminary lemma}

\begin{lem}[\cite{cmf}, Lemma 2.1]\label{lemind} Suppose that $x$
is a singular point of $\Xi$, corresponding to a line bundle $L\in
\PIC^{2g-2}(\CT)$, and there exist $2n$ points
$p_1,\tt(p_1),\ldots,p_n,\tt(p_n)$ of $\CT$ imposing independent
conditions on $H^0(L)$; i.e. if $D=\sum_{i=1}^n(p_i+\tt(p_i))$,
then $h^0(\CT,L(-D))=0$.  Then $\MUL_x \Xi=h^0(\CT,L)/2$.
\end{lem}

\PRF.  With $D$ as above, let $\LL_{D;(1,\ldots,1)}$ be the deformation of $L$
defined in \ref{gooddef}.  By Riemann-Roch, we have
$h^0(L(D))=h^0(\omega_{\CT} \otimes L^{-1}(-D))+2n$. The fact that $L
\in P$ implies that $\omega_{\CT} \otimes L^{-1}\cong \tt^*L$, and
since $D$ is $\tt$-invariant, there is an isomorphism
$H^0(L(-D))\cong H^0(\tt^*L(-D))$ given by $s\mapsto \tt^*s$. 
Hence, $h^0(L(D))=2n$. 

 According to Lemma \ref{Dn1n21}, there is a long exact sequence
\[0\lra H^0(L)\lra H^0(L(D))\stackrel{E}{\lra} H^0(L(D)\otimes
\OO_{D}) \stackrel{\partial_{L;D}}{\lra}H^1(L)\ldots\] and a map
$A_1:H^0(L) \rightarrow H^0(L(D)\otimes \OO_{D})$ such that
$W_1=\ker (\partial_{L;D}\circ A_1)$.  Recall that
$$A_1(s)=(-s(p_1),s(\tau(p_1)),\ldots,-s(p_{n}),
s(\tau(p_{n}))).$$ As observed in the corollary to Lemma
\ref{Dn1n21}, it follows from the above formula that
$\ker(A_1)=H^0(L(-D))=0$. 

On the other hand, since $h^0(L)=h^0(L(D))$, it follows that $\ker
(\partial_{L;D})=0$.  Hence $W_1= \ker (\partial_{L;D}\circ
A_1)=0$, and so $d_k=0$ for all $k\ge 1$.  By Lemma
\ref{lemprymmult} and Proposition \ref{length} (a), $h^0(L)/2 \le
\MUL_x \Xi \le \ell((R^1\pi_* \LL)_{s_0})/2 =d_0/2=h^0(L)/2$. 
Thus $\MUL_x \Xi =h^0(L)/2$. \QED

\subsection{Proof of Theorem \ref{teosv}}

(a) $\iff$ (b).  By the Riemann singularity theorem,
$\MUL_x\THT=h^0(L)$.  Since $\THT \cap P = 2 \cdot \Xi$, it is
clear that $\MUL_x \Xi = (\MUL_x \THT)/2$ if and only if $T_xP
\nsubseteq C_x \THT$. 

\remskip \noindent(b) $\Rightarrow$ (c).  Suppose $\MUL_x \Xi >
h^0(L)/2$. Then by Lemma \ref{lemind}, every choice of $2n$ points
$p_1,\tt(p_1),\ldots,p_n,\tt(p_n)$ of $\CT$ do not impose
independent conditions on $H^0(L)$; i.e. if
$D=\sum_{i=1}^n(p_i+\tt(p_i))$, then $h^0(\CT,L(-D))>0$.  By Lemma
\ref{shpair}, $L=\pi^*M \otimes \OO_{\CT}(B)$,
$h^0(C,M)>h^0(\CT,L)/2$, $B \ge0$, and $B \cap \tau^*B =
\emptyset$. 

\remskip \noindent(c) $\Rightarrow$ (a).  Suppose $L=\pi^*M
\otimes \OO_{\CT}(B)$, $h^0(C,M)>h^0(\CT,L)/2$, $B \ge0$, and $B
\cap \tau^*B = \emptyset$.  Then by Proposition \ref{length} (c),
$\MUL_x \Xi\ge h^0(M)>h^0(L)/2$. 

\remskip \noindent Finally, suppose $L=\pi^*M \otimes
\OO_{\CT}(B)$, $h^0(C,M)>h^0(\CT,L)/2$, $B \ge0$, and $B \cap
\tau^*B = \emptyset$.  Then $M$ and $\OO_{\CT}(B)$ are unique up
to isomorphism by Lemma \ref{uniqueM}, and $h^0(\OO_{\CT}(B))=1$
by Lemma \ref{shpair}. \QED

\section{Proof of Theorem \ref{teo2}}

The basic aim of the proof is to find a deformation of $L$ lying
in the Prym variety, for which $d_1= 2h^0(C,M)-h^0(\CT, L)$, and
$d_2=0$.  One then uses Proposition \ref{length} (c) to conclude.
The computations needed to prove the theorem are quite lengthy,
and consequently, the proof will be broken down into five parts as
follows. 

In Section \ref{1.2} we will fix the class of deformations to be
used in the proof, and establish some preliminary results on
linear systems associated to $L$. In Section \ref{1.3} we will
give a description of the space of sections of $L$ lifting to
first order---a necessary computation for the subsequent sections.
In Section \ref{1.4} we will consider sections lifting to second
order, and we will show that a section lifting to second order
must vanish along a chosen divisor $D$.  In Section \ref{1.5}, we
will show that any section lifting to second order which vanishes
along $D$ must be the zero section, and finally in Section
\ref{1.6} we will complete the proof of the theorem.

\subsection{Preliminaries}\label{1.2}

We will use the following notation.  Let $h^0(\CT ,L)=2n$,
$h^0(C,M)=n_1$ and $h^0(\CT,L)-h^0(C,M)=n_2$.  With regard to
Lemma \ref{length}, we note that $n_1-n_2=2h^0(M)-h^0(L)$. Let
$b\in H^0(\OO_{\CT}(B))$ be a section such that $(b)_0=B$. Let
$D'=\sum_{i=1}^nq_i$, where the $q_i$ are general points of $C$,
and let $\pi^*D'=D= \sum_{i=1}^n(p_i+\tt(p_i))$, where
$\pi^{-1}(q_i)=\{p_i,\tt(p_i)\}$.  With this notation fixed, let
$\LL$ be a family of deformations of $L$ parameterized by $\CX^n$,
whose fiber $\LL_a$ over a point $a\in \CX^n$ is the deformation
$\LL_{D;a}$.

Recall that we are setting $\CT_k=\CT \times \SPEC
\CX[t]/(t^{k+1})$, and denoting by $\LL_k$ the restriction of
$\LL$ to $\CT_k$. We will denote by $W_i(a)$ the image of the map
$H^0(\LL_{a;i})\rightarrow H^0(L)$ induced from the exact sequence
$$0\rightarrow \LL_{a;i-1} \rightarrow \LL_{a;i} \rightarrow L
\rightarrow 0,$$ and by $d_i(a)$ the dimension of this space.  In
other words, for a given $a$, $W_i(a)$ is the space of sections
lifting to order $i$.

Using Lemma \ref{conditions}, and its corollaries we can compute
the dimensions of some pertinent linear systems.  We remark that
for a line bundle $L\in \PIC^{2g-2}(\CT)$ corresponding to a point
$x\in P$, and for a divisor $E$ on $\CT$, the Riemann-Roch theorem
takes the following form: 
$$
h^0(L(E))-h^0(L(-\tt^*E))=\deg(E),
$$
since $\omega_\CT \otimes L^{-1}\cong \tt^*L$, 
and the map $\tt^*:H^0(\tt^*L(-E))\ra H^0(L(-\tt^*E))  $ is an
isomorphism.

\begin{lem}\label{linsys}  In the notation above
\begin{itemize}
\item[\textnormal{(a)}] $h^0(L)=2n$;
\item[\textnormal{(b)}] $h^0(L(-D))=n_1-n$, 
  and $H^0(L(-D))=\pi^*H^0(M(-D'))\cdot b$;
\item[\textnormal{(c)}] $h^0(L(D))=n_1+n$;
\item[\textnormal{(d)}] $h^0(M(-2D'))=
h^0(L(-2D)) =h^0(L(-D-\sum_{i=1}^n p_i))=0$;
\item[\textnormal{(e)}] $h^0(L(2D))=4n$;
\item[\textnormal{(f)}] $h^0(L(2D-\sum_{i=1}^np_i))=3n$;
\item[\textnormal{(g)}] $h^0(L(D+\sum_{i=1}^{n_1-n}\tt(p_i)))=h^0(L(D))$,
so that the natural inclusion induces an isomorphism 
$H^0(L(D))\cong H^0(L(D+\sum_{i=1}^{n_1-n}\tt(p_i)))$;
\item[\textnormal{(h)}]
$h^0(L(-\sum_{i=1}^n\tau(p_i)))=n$;

\item[\textnormal{(i)}] $h^0(L(\sum_{i=1}^np_i))=h^0(L)$,
so that the natural inclusion induces an isomorphism 
$H^0(L) \cong H^0(L(\sum_{i=1}^np_i))$;
\item[\textnormal{(j)}]
$h^0(L(-\sum_{i=1}^n\tau(p_i)-\sum_{i=1}^{n_2} p_i)=h^0(L(-D))$
so that the natural inclusion induces an isomorphism 
$H^0(L(-D))\cong H^0(L(-\sum_{i=1}^n\tau(p_i)-\sum_{i=1}^{n_2} p_i))$.

\end{itemize}
\end{lem}
\PRF.  
(b)  follows from Corollary 
\ref{cor313} (b) with $k=n$.
(c) 
By Riemann-Roch, $h^0(L(D))
-h^0(L(-D))=2n$. 
 Hence, 
$h^0(L(D))=n_1-n+2n=n_1+n$.  
(d) 
follows from Corollary \ref{cor313} (c), 
with $k=k_1=n$.
(e) follows from (d) by Riemann-Roch.
(f)
By Corollary \ref{cor313} (c) with $k=k_1=n$,
$H^0(L(-D-\sum_{i=1}^n p_i))=0$.  (f) then follows
by Riemann-Roch.
(g) 
By Corollary \ref{cor313} (c) with $k=n$ and $
k_1=n_1-n$,
$H^0(L(-D-\sum_{i=1}^{n_1-n}p_i))=0$.
Then by Riemann-Roch,
$h^0(L(D+\sum_{i=1}^{n_1-n}\tt(p_i)))=2n+(n_1-n)=n+n_1=h^0(L(D))$.
(h) The points $\tt(p_1),\ldots,\tt(p_n)$ are general.
(i)  follows from (h)
by Riemann-Roch.
(j) By Corollary \ref{cor313} (b) with $k=n_2$, 
$H^0(L(-\sum_{i=1}^{n_2}(p_i +\tau(p_i))))=
\pi^*H^0(M(-\sum_{i=1}^{n_2}q_i))\cdot b$, and has dimension
$2n-2n_2$.  The  same argument as in 
Corollary \ref{cor313} (c) will show that 
the remaining points 
impose $n-n_2$ conditions.
Hence $h^0(L(-\sum_{i=1}^n\tau(p_i)-\sum_{i=1}^{n_2} p_i))
=2n-2n_2-(n-n_2)=n-n_2=n_1-n=h^0(L(-D))$.
\QED
\remskip

\subsection{Sections lifting to first order}\label{1.3}

\noindent We are now ready to study the sections of $L$ which lift
to first order.  We would like to find some Zariski open subset
$\Omega_1\subseteq \CX^n$ such that for all $a\in \Omega_1$,
$d_1(a)=n_1-n_2$.  In addition, in order to make the second
order computations easier, we will want to understand the relation between
$W_1(a)$ and $H^0(L(-D))$.

\begin{pro}\label{firstorder} 
Let $\Omega_1 =\bigcap_{i=1}^n
\{a_i\ne 0 \} \subseteq \CX^n$.  Then for all $a\in \Omega_1$,
\begin{itemize} 

\item[\textnormal{(a)}] $H^0(L(-D))\subseteq
W_1(a)$; 

\item[\textnormal{(b)}] $d_1(a)=n_1-n_2$;

\item[\textnormal{(c)}] $H^0(L(-\sum_{i=1}^n\tt(p_i)))\cap
W_1(a)=H^0(L(-D))$. 

\end{itemize}
\end{pro}

\PRF.  By Lemma \ref{Dn1n21}, there is a long exact sequence, \[ 0
\lra H^0(L) \lra H^0(L(D)) \stackrel{E_1}{\lra} H^0(L(D)\otimes
\OO_{D}) \stackrel{\partial_{L;D}}{\lra} H^1(L) \ldots\] and a map
$A_1:H^0(L) \rightarrow H^0(L(D)\otimes \OO_{D})$ such that
$W_1(a)=\ker (\partial_{L;D}\circ A_1)$.  Recall that $$
A_1(s)=(-a_1s(p_1),a_1s(\tau(p_1)),\ldots,-a_{n}s(p_{n}),
a_{n}s(\tau(p_{n}))). $$ Thus $\ker(A_1)=H^0(L(-D))$ on
$\Omega_1$, proving (a).  To prove (b), consider the following: 
\begin{eqnarray}\label{kereq}
 \nonumber 
n_1-n_2\ \le \ \dim (W_1(a))&=&\dim(\IM(A_1)\cap \ker (\partial_{L,D}))
+\dim(\ker(A_1))\\ \nonumber
                          &\le   & \dim(\ker (\partial_{L,D}))+\dim(\ker(A_1))\\
                          &=   &(h^0(L(D))-h^0(L))+h^0(L(-D))\\
            \nonumber               & =  &n_1-n_2.
\end{eqnarray}

Finally, to prove (c), consider a section $s\in
H^0(L(-\sum_{i=1}^n\tau(p_i)))\cap W_1(a)$.  $s\in W_1(a)$ is
equivalent to $A_1(s)\in \IM(E_1)$, and it is clear that if $s \in
H^0(L(-\sum_{i=1}^n\tau(p_i)))$, then $A_1(s)\in
E_1(H^0(L(D-\sum_{i=1}^n\tau(p_i))))$. But
$$H^0(L(D-\sum_{i=1}^n\tau(p_i)))=H^0(L(\sum_{i=1}^np_i)),$$ and
we have seen in Lemma \ref{linsys} (i) that
$H^0(L(\sum_{i=1}^np_i))=H^0(L)$.  Therefore $A_1(s)\in
E_1(H^0(L))=0$, and it follows that $s\in H^0(L(-D))$. \QED \remskip

\begin{rem}\label{keyrem} Since the inequalities in (\ref{kereq})
are all equalities, $ \ker (\partial_{L,D})\subseteq \IM(A_1)$. In
fact, $ \ker (\partial_{L,D})=A_1(W_1(a))$, so that
$$ \ker (\partial_{L,D}) = 
\{   
(-a_1s(p_1),a_1s(\tau(p_1)),\ldots,-a_{n}s(p_{n}),
a_{n}s(\tau(p_{n})))|s\in W_1(a)
\}.$$
This will be important in later computations where we will exploit
the fact that $ \ker (\partial_{L,D})$ does not depend on the
$a_i$, whereas $W_1(a)$ does.  \end{rem}

\subsection{Second order lifts: a necessary condition}\label{1.4}

\noindent We are now in a position to consider $W_2$, the space of
sections lifting to second order.  Our eventual goal will be to show 
$W_2=0$. We begin with the following proposition.

\begin{pro}\label{pro531} 
There is a nonempty Zariski open subset
$\Omega_2\subseteq \Omega_1$ such that for all $a \in \Omega_2$,
$W_2(a)\subseteq H^0(L(-D))$.  Therefore, $$W_2(a)\subseteq
H^0(L(-D))\subseteq W_1(a).$$ 
\end{pro}

\PRF.  
By
Lemma \ref{Dn1n22}, 
there is a long exact sequence,
\[
0 \lra H^0(L) \lra H^0(L(2D)) \stackrel{E_2}{\lra} H^0(L(2D)\otimes \OO_{2D}) 
\stackrel{\partial_{L;2D}}{\lra} H^1(L) \ldots\\
\]
and a map $A_2:H^0(\LL_1)
\rightarrow H^0(L(2D)\otimes \OO_{2D})$
such that $\IM(H^0(\LL_2)\rightarrow H^0(\LL_1))=\ker
(\partial_{L;2D}\circ A_2)$.  
In other words, 
a section $s\in H^0(L)$ lifts to second order if and only if there exists
some first order lift $s+\ss^{(1)}t\in H^0(\LL_1)$ such that
$A_2(s+\ss^{(1)}t)\in \IM(E_2)$.  Recall that 
{\small
$$
A_{2}(s+\sigma^{(1)}t)=(0,-a_{1}\sigma^{(1)}(p_{1}),\ a_{1}^{2}s(\tau (p_{1})),
\ a_{1}\sigma^{(1)}(\tau (p_{1})) +a_{1}^{2}\frac{ds}{dz}(\tau 
(p_{1})),\ldots),
$$}
and for $\ffi\in H^0(L(2D))$, 
$$
E_2(\ffi)=(\ffi(p_1),\frac{d\ffi}{dz}(p_1),\ffi(\tau(p_1)),\frac{d\ffi}{dz}
(\tau(p_1)),\ldots),
$$
so that if a section $s$ lifts to second order, then there must be
some section $\ffi \in H^0(L(2D))$ such that $\ffi(p_i)=0$ and
$\ffi(\tt(p_i))=a_i^{2}s(\tau(p_i))$, for all $i$. 

Now let's examine this condition.  Let $F$ be the family of linear
subspaces of $\CX^n$ defined by
$F_a=\{(a_1^2s(\tt(p_1)),\ldots,a_n^2s(\tt(p_n)))\in \CX^n|s\in
W_1\}$, and let $V_2=\{ (\ffi(\tau(p_1)),\ldots,
\ffi(\tau(p_{n})))\in \CX^n|\ffi\in
H^0(L(2D-\sum_{i=i}^{n}p_i))\}$. If a section $s$ lifts to second
order, then $(a_1^2s(\tt(p_1)),\ldots,a_n^2s(\tt(p_n)))\in F_a
\cap V_2$. 

I claim there is a nonempty Zariski open set $\Omega_2\subseteq
\Omega_1$ such that for all $a\in \Omega_2$, $F_a\cap V_2=0$.  It
follows that if a section lifts to second order, it must vanish at
$\tt(p_i)$ for all $i$.  Since $H^0(L(-\sum_{i=1}^n\tt(p_i)))\cap
W_1(a)=H^0(L(-D))$, this means that a section lifting to second
order must be in $H^0(L(-D))$, which completes the proof of the
proposition. \QED

\remskip 
Now we must address the unproven claim. 

\begin{lem}
There is a nonempty Zariski open set $\Omega_2\subseteq \Omega_1$
such that if $a\in \Omega_2$, then $F_a \cap V_2=0$. 
\end{lem}

\PRF.  This will be an application of Proposition \ref{linalg}. We
begin by introducing some notation: with respect to the basis we have
been using for $H^0(L(D)\otimes \OO_D)$, let $pr_2:H^0(L(D)\otimes
\OO_D)\ra \CX^n$ be the projection onto the even factors.  Let
$V_1=pr_2(\ker(\partial_{L,D}))$.

\begin{cla}In the above notation, and for $a\in \Omega_1$,
$\dim(F_a)=n-n_2$.  Furthermore, $F_a =
\{(a_1v_1,\ldots,a_nv_n)\in \CX^n|(v_1,\ldots,v_n)\in V_1\} $.
\end{cla}

\PRF.  Since $H^0(L(-\sum_{i=1}^n\tt(p_i)))\cap
W_1(a)=H^0(L(-D))$, it follows that $F_a\cong W_1(a)/H^0(L(-D))$.
Hence $\dim(F_a)=(n_1-n_2)-(n_1-n)=n-n_2$.  The second statement
is a direct consequence of Remark \ref{keyrem} which implies that
$V_1=pr_2(\ker(\partial_{L,D}))= \{
\left(a_1s(\tau(p_1)),\ldots,a_{n}s(\tau(p_{n}))\right) \ | \ s
\in W_1(a)\}$. \QED

\begin{cla}\label{cla534}

In the above notation, $\dim(V_2)=n_2$.  Furthermore, if
$v=(v_1,\ldots,v_n)\in V_2$, and $v_i=0$ for $i\in I\subseteq
\{1,\ldots,n\}$ with $|I|=n_2$, then $v=0$. \end{cla} \PRF.  $V_2
\cong H^0(L(2D-\sum_{i=1}^{n}p_i))/ H^0(L(D))$. By Lemma
\ref{linsys} (c) and (f), $h^0(L(2D-\sum_{i=1}^{n}p_i))-
h^0(L(D))=3n-(n+n_1)=2n-n_1=n_2$, and hence $\dim(V_2)=n_2$. By
Lemma \ref{linsys} (g), $H^0(L(2D-\sum_{i=1}^{n} p_i -\sum_{i\in
I}\tt(p_{i})))= H^0(L(D+\sum_{i\in I^c}\tt(p_{i})))=H^0(L(D))$,
since $|I^c|=n-n_2=n_1-n$.  Hence if
$v=(\ffi(\tt(p_1)),\ldots,\ffi(\tt(p_n)))\in V_2$, and
$\ffi(\tt(p_i))=0$ for $i\in I$, then $v=0$. \QED \remskip

The proof of the lemma is now just an application of Proposition
\ref{linalg}.  To see this, set $V=\CX^n$. Then
$\dim(V_1)+\dim(V_2)=\dim(V)$, and $F_a =
\{(a_1v_1,\ldots,a_nv_n)\in \CX^n|(v_1,\ldots,v_n)\in V_1\}$, so
that $V$, $V_1$, $V_2$, and $F$ are as in Section 2.2. 
Furthermore, Claim \ref{cla534} shows that $V_2$ satisfies the
conditions of Corollary \ref{cor324}. Hence there exists a
nonempty Zariski open subset $\Omega_2\subseteq \Omega_1$ such
that if $a\in \Omega_2$, then $F_a \cap V_2=0$. \QED

\subsection{Second order lifts: a sufficient condition
}\label{1.5}

Now that we have established that $W_2(a)\subseteq
H^0(L(-D))\subseteq W_1(a)$, we can focus our attention on
sections in $H^0(L(-D))$.  This is a great advantage, as we know
the exact form of the first order lifts of such sections.  In
order to take full advantage of this information, we prove the
following lemma which addresses an important special case.  First,
let us define the following notation. Let
$\{e_1^*,\ldots,e_{4n}^*\}$ be the dual basis to the basis we have
been using for $H^0(L(2D)\otimes \OO_{2D})$.  Let $H = \{v\in
H^0(L(2D)\otimes \OO_{2D})|e_{i+4j}^*(v)=0 \textnormal{ for }
i=1,3,4 \textnormal{ and } 0\le j \le n-1\}$;  i.e.
$H=\{(0,*,0,0,\ldots,0,*,0,0) \in H^0(L(2D)\otimes \OO_{2D})\}$.

\begin{lem} 
If $s +\ss^{(1)}t\in \IM(H^0(\LL_2) \ra H^0(\LL_1))$ and
$A_2(s+\ss^{(1)}t)\in H$, then $A_2(s+\ss^{(1)}t)=0$. 
\end{lem}

\PRF.  If $s+\ss^{(1)}t$ lifts to second order, then
$A_2(s+\ss^{(1)}t)=E_2(\ffi)$ for some $\ffi\in H^0(L(2D))$. Due
to the form of $H$, we can see that $\ffi \in
H^0(L(2D-D-\sum_{i=1}^n\tt(p_i)))= H^0(L(\sum_{i=1}^np_i))$. By
Lemma \ref{linsys} (i), $H^0(L(\sum_{i=1}^np_i))=H^0(L)$, so that
$\ffi \in H^0(L)$, and consequently $E_2(\ffi)=0$. Thus
$A_2(s+\ss^{(1)}t)=0$.\QED 
\remskip

With this we will prove the next proposition.

\begin{pro}\label{secondorder}
There is a Zariski open subset $\Omega_3\subseteq \Omega_2$ such
that for all $a \in \Omega_3$, $W_2(a)=0$. 
\end{pro}

\PRF.  Let $s\in W_2(a)$.  By Proposition \ref{pro531}, $s\in
H^0(L(-D))$.  Let $s+\ss^{(1)}t$ be the standard lifting given in
Lemma \ref{canon}, and recall that on an open set $U_i \subseteq
\CT$ in our cover,

\[
\ss_i^{(1)}= \left\{ \begin{array}{ll}
a_is/z, &\textnormal{if } p_i\in U_i;\\
-a_is/z, &\textnormal{if } \tau(p_i)\in U_i;\\
0,   &\textnormal{otherwise}.
\end{array}\right. \]
A general lifting of $s$ will be given by $s+(\ss^{(1)}+\ffi)t$
for some $\ffi \in H^0(L)$.  Observe that
\begin{equation}\label{eq543} A_{2}(s+(\ss^{(1)}+\ffi)t)=
(0,-a_{1}^2(s/z)(p_{1})-a_1\ffi(p_1), \ 0\ ,
a_1\ffi(\tt(p_{1})),\ldots), \end{equation} since
$(s/z)(\tt(p_i))=(ds/dz)(\tt(p_i))$ when $s(\tt(p_i))=0$. 

I claim there is a $\psi \in W_1$ such that
$\psi(\tt(p_i))=\ffi(\tt(p_i))$ for all $i$.  Indeed if
$A_{2}(s+(\ss^{(1)}+\ffi)t) \in \IM(E_2)$, then it must be in
$E_2(H^0(L(D)))$, since $s\in H^0(L(-D))$. Consider the following diagram,
\[
\begin{array}{ccc}
H^0(L(2D))& \stackrel{E_2}{\lra}&H^0(L(2D)/L)\\
\bigcup  &                     & \bigcup\\
H^0(L(D))& \stackrel{E_1}{\lra}&H^0(L(D)/L),\\
\end{array}
\]
where the inclusion on the right, in our chosen bases, is given by
$$
(x_1,x_{\tt(1)},\ldots,x_n,x_{\tt(n)})
\mapsto (0,x_1,0,x_{\tt(1)},\ldots,0,x_n,0,x_{\tt(n)}).
$$
It follows that 
$A_2(s+\ss^{(1)}t) \in H^0(L(D)/L) \subseteq H^0(L(2D)/L)$, 
so that
$$(-a_{1}^2(s/z)(p_{1})-a_1\ffi(p_1),
a_1\ffi(\tt(p_{1})),\ldots)
\in E_1(H^0(L(D))).$$
Recall from Remark \ref{keyrem} that $\IM(E_1)=A_1(W_1)$, so that
there must be a $\psi\in W_1$ such that $A_1(\psi)=
(-a_1\psi(p_1),a_1 \psi(\tt(p_1)), \ldots)=
(-a_{1}^2\frac{s}{z}(p_{1})-a_1\ffi(p_1),
a_1\ffi(\tt(p_{1})),\ldots)$. Hence,
$\psi(\tt(p_i))=\ffi(\tt(p_i))$ for all $i$, establishing the
claim. 

Since $\psi$ lifts to first order, $\psi t$ lifts to second order,
so that $s+(\ss^{(1)}+\ffi)t -\psi t$ also lifts to second order. 
But then $A_2(s+(\ss^{(1)}+\ffi)t -\psi t) \in H$, so that the
above lemma implies that $s+(\ss^{(1)}+\ffi)t -\psi t \in \ker
(A_2)$.  Setting $\rho = \psi-\ffi$, we have that $A_2(s
+(\ss^{(1)} + \rho)t)=0$.  In other words, using equation
(\ref{eq543}), if $s$ lifts to second order, then there is a
section $\rho \in H^0(L)$ such that $\rho(\tt(p_i))=0 $ and
$\rho(p_i)=a_i\frac{s}{z}(p_i)$ for all $i$.

Now let's examine this condition.  Let $F$ be the family of linear
subspaces of $\CX^n$ defined by $F_a=\{ (a_{1}
\big(\frac{s}{z}\big)(p_{1}),\ldots,a_{n}
\big(\frac{s}{z}\big)(p_{n})) \in \CX^n|s\in H^0(L(-D))\}$, and
let $V_2=\{ (\rho(p_1),\ldots, \rho(p_{n})) \in \CX^n|\rho\in
H^0(L(-\sum_{i=1}^{n}\tt(p_i)))\}$. If a section $s$ lifts to
second order, then $ (a_{1}
\big(\frac{s}{z}\big)(p_{1}),\ldots,a_{n}
\big(\frac{s}{z}\big)(p_{n})) \in F_a \cap V_2$. 

I claim there is a nonempty Zariski open subset $\Omega_3\subseteq
\Omega_2$ such that for all $a\in \Omega_3$, $F_a\cap V_2=0$.  It
follows that if a section lifts to second order, it must vanish to
second order at $p_i$ for all $i$, so that $s \in
H^0(L(-D-\sum_{i=1}^n p_i)).$ By Lemma \ref{linsys} (d),
$H^0(L(-D-\sum_{i=1}^n p_i))=0$, and hence $s=0$. \QED

\remskip
Now we must  address the unproven claim.
\begin{lem}
There is a nonempty Zariski open subset $\Omega_3\subseteq
\Omega_2$ such that if $a\in \Omega_3$, then $F_a\cap V_2=0$. 
\end{lem}

\PRF.  This will be an application of Proposition \ref{linalg}. To
begin, let $V_1= \{ (\frac{s}{z}(p_1),\ldots,\frac{s}{z}(p_n)) \in
\CX^n| s\in H^0(L(-D))\}$, so that $V_1=F_{(1,\ldots,1)}$. It is
clear that $V_1 \cong H^0(L(-D))/ H^0(L(-D-\sum_{i=1}^n
p_i))=H^0(L(-D))$, since
 $H^0(L(-D-\sum_{i=1}^n p_i))=0$.

\begin{cla}
In the above notation and for $a\in \Omega_1$, 
$
\dim(F_a)=n_1-n=n-n_2$.  
\end{cla}

\PRF.  On $\Omega_1$, $F_a \cong V_1 \cong H^0(L(-D))$. Thus,
$\dim(F_a)=h^0(L(-D))=n_1-n$. \QED

\begin{cla}\label{cla545}

In the above notation, $\dim(V_2)=n_2$.  Furthermore, if
$v=(v_1,\ldots,v_n)\in V_2$, and $v_i=0$ for $i\in I\subseteq
\{1,\ldots,n\}$ with $|I|=n_2$, then $v=0$. \end{cla}

\PRF.  $V_2 \cong H^0(L(-\sum_{i=1}^n\tau(p_i))/H^0(L(-D))$. By
Lemma \ref{linsys} (h) and (b), $h^0(L(-\sum_{i=1}^n\tau(p_i)))
=n$ and $h^0(L(-D))=n_1-n$, so that
$\dim(V_2)=n-(n_1-n)=2n-n_1=n_2$.  By Lemma \ref{linsys} (j),
$H^0(L(-\sum_{i=1}^n\tau(p_i)-\sum_{i\in I} p_i))=H^0(L(-D))$.
Hence if $v=(\ffi(p_1),\ldots,\ffi(p_n))\in V_2$, and
$\ffi(p_i)=0$ for $i\in I$, then $v=0$. \QED 
\remskip

The proof of the lemma is now just an application of Proposition
\ref{linalg}.  To see this, set $V=\CX^n$. Then
$\dim(V_1)+\dim(V_2)=\dim(V)$, and $F_a =
\{(a_1v_1,\ldots,a_nv_n)\in \CX^n|(v_1,\ldots,v_n)\in V_1\}$, so
that $V$, $V_1$, $V_2$, and $F$ are as in Section 2.2. 
Furthermore, Claim \ref{cla545} shows that $V_2$ satisfies the
conditions of Corollary \ref{cor324}. Hence there exists a
nonempty Zariski open set $\Omega_2\subseteq \Omega_1$ such that
if $a\in \Omega_2$, then $F_a \cap V_2=0$. \QED

\subsection{Proof of  Theorem \ref{teo2}}\label{1.6}

Let $a \in \Omega_3$, and consider the deformation $\LL_{D;a}$. 
It follows from Proposition \ref{firstorder} that
$d_1=n_1-n_2=2h^0(M)-h^0(L)$.  Proposition \ref{secondorder}
implies that $W_2=0$, and hence that $d_2=0$.  Due to Lemma
\ref{length} (c), $\MUL_x\Xi=h^0(M)$.\QED

\section{Consequences of Theorem \ref{teo2}}

\subsection{ Upper bounds on multiplicity of singularities}

For a ppav $(A,\Theta)$, 
let $\singd_k\Theta=\{x\in \sing\Theta\ |\ \MUL_x \Theta\ge k\}$. 
A result of Koll\'ar's \cite{kbound} shows that if
$\dim(A)=d$, then $\dim(\singd_k\Theta)\le
d-k$.  Generalizing a result of Smith and Varley \cite{svbound}, Ein and Lazarsfeld \cite{elbound} showed  that 
$\dim(\singd_k\Theta)=d-k$ only if $(A,\Theta)$ splits as a 
$k$-fold product.  Thus for an irreducible Prym variety  associated to a connected \'etale double cover of a smooth curve $C$ of genus $g$, and 
$x\in \sing \Xi$,  $\MUL_x \Xi\le g-2=\dim(P)-1$.  
Using Theorem \ref{teo2}, we will improve these estimates
for Prym varieties.  To begin,
we will prove the following lemma on double covers of hyperelliptic curves: 

\begin{lem}[\cite{cmf}, Lemma 3.5]\label{lemhype}
Let $\pi:\CT \ra C$ be a connected \'etale double cover of a
smooth curve $C$. If $\CT$ is hyperelliptic, then $C$ is
hyperelliptic.  Furthermore, if $\tilde{G}$ is the line bundle
corresponding to the $g^1_2$ on $\CT$, and $G$ is the line bundle
corresponding to the $g^1_2$ on $C$, then $\textnormal{Norm}(\tilde{G})\cong G$
and $\pi^*G\cong \tilde{G}^{\otimes 2}$. 
\end{lem}

\PRF.  Suppose $p_1+p_2$ and $p_3+p_4$ are general 
in $|\tilde{G}|$.  Since Norm preserves linear equivalence, 
$\pi(p_1)+\pi(p_2)\sim \pi(p_3)+\pi(p_4)$.  Thus $C$ is hyperelliptic, and since there is a unique $g^1_2$ on $C$, 
$\textnormal{Norm}(\tilde{G})\cong G$.
Now 
let $f:\CT\ra \PR$ be the morphism corresponding to
the $g^1_2$ on $\CT$.  Then
$f\circ \tau$ is also a finite degree $2$ morphism of
$\CT$ to $\PR$, and since there is a unique $g^1_2$ on $\CT$, this
implies that $\tt(p_1)+\tt(p_2) \sim p_1+p_2$.  Thus $\pi^*G
\cong 
\pi^* \textnormal{Norm}\tilde{G}\cong 
\OO_{\CT}(p_1+p_2+\tt(p_1)+\tt(p_2))\cong \tilde{G}^{\otimes
2}$. \QED \remskip

Recall the definition of the Clifford index:
$$\cliff(C)=\min\{\deg(D)-2\dim|D| :  h^0(D)\ge 2, h^1(D)\ge
2\}.$$ If $x\in \sing \Xi$, then $h^0(L)\le
\deg(L)/2-\cliff(\CT)/2+1$, and $h^0(M)\le
\deg(M)/2-\cliff(C)/2+1$. Letting $W^r_d(C)$ denote the variety of line bundles 
$L$ on $C$ such that $\deg(L)=d$ and $h^0(L)>r$,
Martens' theorem \cite{martens} states that
if $2\le d \le g-1$ and $0<2r\le d$, then
$\dim(W^r_d(C))\le d-2r$, with equality holding only if
$C$ is hyperelliptic.  Since 
$\cliff(C)\le d-2r$, these inequalities yield essentially the same information.

\begin{cor}\label{upbound}
If $x\in \sing\Xi$, then $\MUL_x \Xi \le (g+1)/2$.
If $P$ is irreducible, then
$\MUL_x \Xi \le g/2=(\dim(P)+1)/2$.
More precisely, 
 suppose $g\ge 5$, and 
 let $Z$ be an irreducible component of
$\singd_k\Xi$.  For $x\in Z$, let $L_x$ be the corresponding line bundle.
\begin{enumerate}
\item[\textnormal{(a)}]
Suppose for a general $x\in Z$, 
$\MUL_x\Xi=h^0(L_x)/2$, i.e. $T_x P\nsubseteq C_x\THT$. Then $k\le g/2-\cliff(\CT)/4$,
 and $\dim(Z)\le 2g-4k$.
 If $\dim(Z)= 2g-4k$, then
 $C$ is hyperelliptic.
 
\item[\textnormal{(b)}]
 Suppose for a general $x\in Z$, 
$\MUL_x\Xi>h^0(L_x)/2$, i.e. $T_x P\subseteq C_x\THT$.  Then $k\le (g+1)/2-\cliff(C)/2$, and
$\dim(Z)\le g-2k+1$.
If $\dim(Z)= g-2k+1$, then $C$ is 
hyperelliptic.  If we suppose moreover that
$C$ is not hyperelliptic, 
and $\dim(Z)>0$, then $\dim(Z)\le g-2k-1$, 
and if 
$\dim(Z) =  g-2k-1$, then 
$C$ is either trigonal, bi-elliptic, or a smooth plane quintic.
\end{enumerate}

\end{cor}

\PRF.  The first statement of the corollary follows immediately from (a) and (b).  The statement for
irreducible Prym varieties then follows from 
Mumford's result \cite{mum1} that if $C$ is hyperelliptic, then 
the Prym variety associated to the double cover is a hyperelliptic Jacobian, or the product of two such Jacobians. 

(a) We have that 
$\MUL_x\Xi=h^0(L)/2$.  Since $\deg(L)=2g-2=g(\CT)-1$, 
by Riemann-Roch, $h^1(L)=h^0(L)\ge 2$.  Thus by Clifford's theorem
$\MUL_x\Xi\le (2g-2)/4 -\cliff(\CT)/4+1/2=g/2-\cliff(\CT)/4$.  In addition, 
we must have that
$Z\subseteq W^{2k-1}_{2g-2}(\CT)$, and so it
follows immediately from Martens' Theorem that 
$\dim(Z)\le 2g-2-2(2k-1)=2g-4k$, with equality holding
only if $\CT$, and hence $C$, is hyperelliptic.

(b)  Now suppose $T_x P \subseteq C_x \THT$, so that $L=\pi^*M
\otimes \OO_{\CT}(B)$, $h^0(C,M)>h^0(\CT,L)/2$, $B \ge0$, $B \cap
\tau^*B = \emptyset$, and $\MUL_x \Xi = h^0(C,M)$.  Since $M$ is special,
Clifford's theorem implies
$h^0(C,M)\le \deg(M)/2 -\cliff(C)/2+1\le (g+1)/2
-\cliff(C)/2$.  Applying Martens' Theorem, 
we see that $\dim(Z)\le (g-1)-2(k-1)= g-2k+1$, with 
equality holding only if $C$ is hyperelliptic.
If we suppose that $\dim(Z)>0$, then we can assume that
$\deg(M)<g-1$, since there are only a finite number of theta characteristics.  
If we assume further that $C$ is not hyperelliptic, then it follows from Mumford's refinement \cite{mum1} of Martens' theorem that
$\dim(Z)\le g-2k-1$, with equality holding only if
$C$ is either trigonal, bi-elliptic, or a smooth plane quintic.
\QED

\remskip
\begin{rem}  The statements regarding the dimension of 
$\singd_k\Xi$ in the corollary above were pointed out by the referee.  It should also be noted that 
in \cite{mum1}, Mumford 
studied a skew symmetric bilinear
pairing 
$$\beta :H^0(L)\wedge H^0(L)\rightarrow H^0(\omega_C\otimes \eta),$$
 and showed by a dimension count
that if $g\ge 5$, and $\dim(\sing\Xi)\ge g-5$, then $\beta$ has a two dimensional isotropic subspace $V\subseteq H^0(L)$.
Such an isotropic subspace gives rise to an 
isomorphism
$L\cong \pi^*M\otimes \OO_\CT(B)$, where 
$h^0(M)\ge 2$, $B\ge 0$, $B\cap \tt^*B = \emptyset$.
As in the proof of (b) above, Mumford 
concluded that if $C$ is not hyperelliptic, then
$\dim(\sing\Xi)\le g-5$.
Generalizing Mumford's work, 
Smith and Varley \cite{sv2} have shown that there exists an 
isotropic subspace of dimension $k$ for $\beta$ if and
only if there exists such a decomposition of $L$ with 
$h^0(M) \ge k$. Thus it seems possible that through further
analysis of the pairing $\beta$, one may be able to improve the bound 
on $\dim(\singd_k\Xi)$.
\end{rem}
\remskip
\begin{rem} The referee has raised the question of whether Ein and Lazarsfeld's bound on the dimension of
$\singd_k\Theta$ for irreducible ppavs is sharp.  I.e. do
there exist irreducible ppavs with  
$\dim(\singd_k\Theta)=d-k-1$?  As an example, it would be interesting to know if there exist irreducible ppavs of dimension five with a point of order four on their theta divisor.
It appears that the techniques of this paper may extend to Prym varieties associated to double covers of stable curves, and hence in the case of an irreducible ppav of dimension  less than or equal to five, it may be possible to answer this question and give a sharp bound on 
$\dim(\singd_k\Theta)$.  This is work in progress.
\end{rem}
\remskip
\begin{rem}
For the Jacobian of a curve, Martens' theorem implies that 
$\textnormal{codim} (\singd_k\Theta) = 2k-1 $
only if the curve is hyperelliptic.
  It is a result of 
Beauville \cite{b1} that
if $(A,\Theta)$ is an irreducible generalized Prym variety, and $\dim (\sing\Xi) \ge g-4=\dim(A)-3$, then 
$(A,\Theta)$ is a hyperelliptic Jacobian
(c.f. Mumford \cite{mum1}).  
Thus at least in dimension 
less than or equal to five,
any irreducible ppav whose theta divisor has  
double points in codimension three is a hyperelliptic Jacobian. In regards to these results,
and part (b) of the Corollary,
the referee has asked to what extent $k$-fold points in codimension $2k-1$ characterize hyperelliptic Jacobians among all irreducible ppavs.
Despite the examples cited, in general for $k>2$ it appears there may be other components in this locus, since
at least in dimension five, the theta divisor of the intermediate Jacobian of a smooth cubic threefold has a triple point, but such a ppav is not a Jacobian (Clemens and Griffiths \cite{cg}).
\end{rem}

\remskip
Recall that given distinct points $p_1,\tt(p_1),\ldots,p_k, \tt(p_k)$ of
$\CT$, we define $D_k=\sum_{i=1}^k(p_i+\tt(p_i))$. Due to Remark
\ref{remcomp}, we have the following upper bound on the multiplicity of a point on the Prym theta divisor: 

\begin{cor}
Suppose $x\in \sing\Xi$ corresponds to the line bundle $L \in
\PIC^{2g-2}(\CT)$.  If there exist $2k$ distinct points
$p_1,\tt(p_1),\ldots,p_k, \tt(p_k)$ of $\CT$ such that
$h^0(L(-D_k))=h^0(L)-2k$, then $$ \MUL_x \Xi \le h^0(L)-k, $$ with
equality holding if and only if $k$ is the largest number with
this property. \QED 
\end{cor}

\subsection{Prym varieties of dimension five}

A Prym variety of dimension five is associated to a double cover
of a genus six curve. For a point $x\in \sing \Xi$, Corollary
\ref{upbound} implies that $\MUL_x\Xi \le 3$; in this section we
will examine exactly which Prym varieties of dimension five have
singular theta divisors with triple points. Theorem
\ref{teogen6} is a direct consequence of the following theorem. 

\begin{teo}\label{teopq}

Suppose $\dim P=5$. If $\singd_3 \Xi \ne \emptyset$, then one of
the following must hold: 

\begin{itemize}
\item[\textnormal{(a)}]
$C$ is a plane quintic and 
$h^0(\OO_{\PRR}(1)|_C \otimes \eta)=1$.  In this case
$\singd_3\Xi=\sing \Xi=\{x\}$ 
consists of a unique point corresponding to the line
bundle
$\pi^*(\OO_{\PRR}(1)|_C)$.  Moreover, $X=C_x \Xi$ is a smooth 
cubic threefold, and 
$(P,\Xi)\cong (JX, \Theta)$.
 
\item[\textnormal{(b)}] $C$ is hyperelliptic,
and either
\begin{itemize}
\item[\textnormal{(i)}]$\CT$ is
hyperelliptic and $(P, \Xi)\cong JC'$ for some hyperelliptic curve
$C'$.  In this case $\dim(\sing\Xi)=2$, and
$\singd_3\Xi=\{x\}$ consists of a unique point corresponding to the line
bundle
$5g^1_2$ on $\CT$;

\item[\textnormal{(ii)}]
$\CT$ is not hyperelliptic, and
$(P,\Xi)\cong JC'\times JC''$, for some hyperelliptic curves
$C'$ and $C''$.  In this case $\dim(\sing\Xi)=3$, and
$\dim (\singd_3\Xi)=1$.
\end{itemize}

\end{itemize}
\end{teo}

\PRF.  (a) Suppose $C$ is not hyperelliptic, and 
$L$ is a line bundle corresponding to a singular
point $x$ of multiplicity $3$.  Since $\deg(L)=10$, by
Clifford's theorem, $h^0(L)\le 6$, with equality holding
only if $\CT$, and hence $C$, is hyperelliptic.  
By Theorems \ref{teosv} and \ref{teo2}, if $h^0(L)=2$, then 
$\MUL_x\Xi\le 2$.  
Thus we may assume $h^0(L)=4$, and by
Theorem \ref{teosv}, 
$L=\pi^*M\otimes \OO_{\CT}(B)$, with $h^0(M)\ge3$, $B\ge 0$
and $B\cap \tt^*B = \emptyset$. By Clifford's theorem applied to $M$, either $C$ is hyperelliptic, or $\deg(M)=5$,
$h^0(M)=3$, and
$L=\pi^*M$.
We must have that $|M|$ is base point free, or else there would be a
$g^2_4$ on $C$.  Thus $|M|$ defines a morphism to $\PRR$, which is
birational since $M$ has prime degree, and is an embedding since
the genus of a smooth plane quintic is $6$.  Hence, $M\cong
\OO_{\PRR}(1)|_C$, and $h^0(\OO_{\PRR}(1)|_C\otimes \eta)=1$,
since $h^0(L)=h^0(M)+h^0(M\otimes \eta)$. A plane quintic has a
unique $g^2_5$ (see \cite{acgh} p.209), and so $L$ can be the only
triple point on $\Xi$. 

Given a smooth plane quintic such that $h^0(\OO_{\PRR}(1)|_C
\otimes \eta)=1$, an elementary argument (c.f. Beauville
\cite{b0}) will show that $\Xi$ has a unique singular point.
Finally, Beauville \cite{b2} (see also Donagi and Smith \cite{ds},
as well as \cite{cmf}), has shown that $X=C_x\Xi$ is a smooth
cubic threefold, and $(P,\Xi)\cong (JX, \Theta)$ as principally
polarized abelian varieties. 

(b) part (i).  If $\CT$ is hyperelliptic, then the proof of
Mumford's theorem \cite{mum1} p.344 implies that $(P,\Xi)\cong
(JC',\Theta')$, for some hyperelliptic curve $C'$ of genus $5$. 
Hence, $\dim(\sing \Theta')=2$, and $\singd_3\Theta'=\{x\}$, where
$x$ corresponds to the unique $2g^1_2$ on $C'$.  On the other
hand, $h^0(\CT, 5g^1_2)=6$, and $\pi_*(5g^1_2)=5g^1_2=\omega_C$,
so that $5g^1_2$ corresponds to a triple point on $\Xi$. 

(b) part (ii).  If $\CT$ is not hyperelliptic, then the proof of
Mumford's theorem implies that $(P,\Xi)\cong (JC'\times JC'',
JC'\times \Theta''+\Theta'\times JC'')$, for some hyperelliptic
curves $C'$ and $C''$. The possible genera for $C'$ and $C''$ are
$1$ and $4$, or $2$ and $3$, respectively. 

In the former case, $\sing \Xi=(\Theta'\times \Theta'') \cup
(JC'\times \sing \Theta'')$, and it follows that
$\dim(\sing\Xi)=\dim (\Theta'')= 3$.  $\singd_3\Xi=\Theta'\times
\singd_2\Theta''$, and $\singd_2\Theta''=\{g^1_2+p\}$, which has
dimension one, so $\dim(\singd_3\Xi)=1$. 

In the latter case, $\sing \Xi=(\Theta'\times \Theta'') \cup
(JC'\times \sing \Theta'')$, and it follows that
$\dim(\sing\Xi)=\dim (\Theta') + \dim (\Theta'')= 3$. 
$\singd_3\Xi=\Theta'\times \singd_2\Theta''$, and thus
$\dim(\singd_3\Xi)=\dim(\Theta')+\dim(\singd_2\Theta'') =1$.\QED

\remskip 
\begin{rem}
The proof above includes a simplification suggested by the referee, who also observed that this theorem is 
deducible from the results of Friedman and the author in \cite{cmf}.  To be precise, the proof 
of Theorem \ref{teopq} only uses the special case of Theorem \ref{teo2} that  
$h^0(L)=2$.  This special case follows from \cite{cmf}
Theorem 2.5,  p.306.
\end{rem}

\remskip 
\begin{rem} Some of the statements in part (b) can be
proven without using the fact that the Prym of a hyperelliptic
curve is a hyperelliptic Jacobian.  Namely, one can show that if
$\CT$ is hyperelliptic, then there is a unique triple point of
$\Xi$, and if $\CT$ is not hyperelliptic, then
$\dim(\singd_3\Xi)\le 1$. Indeed, as observed in the proof above,
if $\CT$ is hyperelliptic, then there is a unique $g^5_{10}$ on
$\CT$, namely $5g^1_2$.  Furthermore,
$\pi_*5g^1_2=5g^1_2=\omega_C$, so that in fact $5g^1_2\in
\singd_3\Xi$.  I claim that there are no triple points with
$h^0(L)=4$.  In fact, since $h^0(M)=3$, we must have $M\ge
2g^1_2$.  But then $\pi^*M\ge 4\tilde{g}^1_2$, so that
$h^0(\pi^*M)\ge 5 >4=h^0(L)$, a contradiction. On the other hand,
if $C$ is hyperelliptic, then $h^0(L)=4$.  Since $h^0(M)=3$, $M\ge
2g^1_2$.  Let $M'=2g^1_2$, so that $L=\pi^*M'\otimes
\OO_{\CT}(B')$, where now we only require that $B'>0$ and
$\deg(B')=2$.  It follows that $\pi_*L=\pi_*\pi^*g^1_2\otimes
\pi_*B'= 4g^1_2\otimes \pi_*B'$, so that if $\pi_*L=\omega_C$,
then $B'$ must lie above the $g^1_2$ on $C$.  Since the $g^1_2$
has dimension one, and there are four choices of $B'$ above each
pair of points in the $g^1_2$, the dimension of triple points can
be at most one. 

\end{rem}

\section{The Prym Canonical Map}
 
Let $\Psi_{\eta}:C \ra \PROJ^{g-2}$ be the Prym canonical map;
i.e. the map induced by the linear system $|\omega_{C}\otimes \eta
|$.  One can easily check that $|\omega_C \otimes \eta|$ has a
base point if and only if $\CT$ is hyperelliptic, and consequently
in this section we will restrict our attention to the case that
$\CT$ is not hyperelliptic. Under this assumption, we will
establish a connection between the tangent cone to a singular
point of $\Xi$ and the Prym canonical image of $C$.  Our eventual
goal will be to determine whether a $k$-dimensional secant variety
to $C$ is contained in the tangent cone to a singular point of
$\Xi$. 

\subsection{Preliminaries on Prym images}

Let $q\in C$, and set $D=p+\tt(p)=\pi^*(q)$. Consider the exact
sequence
$$
0\lra \OO_{\CT} \lra \OO_{\CT}(D) \lra 
\OO_D
\lra 0,
$$
and let $\partial_D$ be the boundary map of the associated long
exact sequence.  Since $\CT$ is not hyperelliptic,
$h^0(\OO_{\CT})=h^0(\OO_{\CT}(p+\tt(p)))=1$, and hence
$\partial_D$ is injective.

 \begin{lem}\label{genprym}
$\Psi_{\eta}(q)= \{ \partial_D(-a,a)\ | \ a \in \CX\}\in \PROJ
(H^1(\OO_C)^-)$.
 \end{lem}
\remskip

\PRF.  Serre duality gives an isomorphism $H^1(\OO_{\CT})\ra
(H^0(\omega_{\CT}))^*$, given by $\aa \mapsto (\omega
\mapsto\sum_{p'\in \CT} \RES_{p'}(\aa\omega))$. From this
description, it is easy to see that this isomorphism induces an
isomorphism $H^1(\OO_{\CT})^-\ra (H^0(\omega_{\CT})^- )^*$, and
hence an isomorphism $\PROJ(H^1(\OO_{\CT})^-)\ra
\PROJ((H^0(\omega_{\CT})^- )^*)$. The map $C\ra \PROJ
(H^1(\OO_{\CT})^-)$ given by $q \mapsto \langle\xi \rangle =\{
\partial_{D}(-a,a)\ | \ a \in \CX\}$ is well defined since
$\partial_D$ is injective. Composing with the duality map gives a
map $\psi:C\ra \PROJ (H^0(\omega_{\CT})^- )^*)$ which by
definition is given by $q \mapsto \langle \omega \mapsto
\sum_{p'\in \CT} \RES_{p'}(\xi\omega)\rangle$.

Letting $\xi = \partial_{\pi^*(q)}(-a,a)$, a computation in the
$\CH$ech complex shows that in the open cover we have been using,
$\xi$ is given as follows: 
\begin{displaymath}
    \xi= \left\{ 
    \begin{array}{ll}
        a/z, & \textnormal{on } U_{0i},\\
        -a/z,&\textnormal{on } U_{\tau(0)j},\\
        0,   & \textnormal{otherwise }.
        \end{array}
        \right.
\end{displaymath} 
Hence, $\sum_{p'\in \CT} \RES_{p'}(\xi\omega)=2a\omega(p)$. It
follows that $q \stackrel{\psi}{\mapsto} \langle \omega \mapsto
\omega(p) \rangle$, which is the definition of $\Psi_{\eta}$. \QED
\remskip

Now consider a point $x\in \sing \Xi$, corresponding to the line
bundle $L\in \PIC^{2g-2}(\CT)$, and consider the deformation
$\LL_{D,a}$ for some $a\in \CX$.  As before, let the transition
functions of $\LL_1$ be denoted by
$\ll_{ij}(t)=\ll_{ij}(1+\aa_{ij}^{(1)}t)$, and let $f:S \ra P$ be the
associated morphism.  Then $f_*:T_{s_0}S \ra T_x P =
H^1(\OO_{\CT})^-$ has image equal to the linear span of
$\aa^{(1)}\in H^1(\OO_{\CT})^-$.  From our computation of
$\aa^{(1)}$ in Section 1, and the description of
$\partial_D(-a,a)$ given in the proof above, we have the
following: 

\begin{lem}
In the above notation, $$f_*(T_{s_0}S) =\langle \aa^{(1)}\rangle
=\langle \partial_D(-a,a)\rangle = \Psi_\eta(q) \in
\PROJ(H^1(\OO_C)^-).$$ \QED 
\end{lem} \remskip

More generally, we can ask for the relation between the Prym
canonical map, and a deformation $\LL_{D;a}$, where $D$ has higher
degree. Let $q_1, \ldots q_k \in C$, and set
$D=\pi^*(\sum_{i=1}^kq_i)$. Let $\LL$ be a family of line bundles
over $C$, parameterized by $\CX^k$, such that $\LL_a$ is a line
bundle associated to $\LL_{D,a}$, and consider the induced family
of maps $f_a :S_a \ra P$.  For all $a$, let $s_0 \in S_a$ be such
that $f_a(s_0)=x$.  Finally, let
$\langle\Psi_{\eta}(q_{1}),\ldots,\Psi_{\eta}(q_{k})\rangle
\subseteq \PROJ^{g-2}$ be the span of the points $\Psi_\eta
(q_{1}),\ldots,\Psi_\eta(q_{k})$.  Extending the proofs of the
first two lemmas by linearity, we immediately have the following:

\begin{lem}
In the above notation,
$$
\langle\Psi_{\eta}(q_{1}),\ldots,\Psi_{\eta}(q_{k})\rangle
=\{ \partial_D(-a_{1},a_{1},\ldots,-a_{k},a_{k})\ | \ 
a 
\in \CX^k \}
\subseteq \PROJ(H^1(\OO_{\CT})^-),$$
and for each $a\in \CX^k$,
$$
(f_{a})_*(T_{s_0}S_a)=\langle \aa_a^{(1)} \rangle=
\langle \partial_D(-a_{1},a_{1},\ldots,-a_{k},a_{k})\rangle \in
\PROJ (H^1(\OO_C)^-).
$$

\QED
 \end{lem}

As a consequence, we have the following proposition:

\begin{pro}\label{teoprymim}Suppose
$\Psi_{\eta}(q_1), \dots ,\Psi_{\eta}(q_k)$ are contained in a unique $(k-1)$-plane.
 Then $\partial_D$ induces a linear inclusion
$\partial_D: \PROJ^{k-1} \ra \PROJ^{g-1}$.
Moreover, if $\MUL_x \Xi =\mu $, then $\partial_D$ gives
 a bijection of sets
\[
\{a \in \PROJ^{k-1}\ | \ (\MUL_{s_0}\Theta_{S_a})/2 >\mu\}
\leftrightarrow
\langle\Psi_{\eta}(q_{1}),\ldots,\Psi_{\eta}(q_{k})\rangle \cap C_x \Xi.
\]
\end{pro}

\PRF.  $\MUL_{s_0}\Theta_{S_a}/2>\mu$ if and only if $\langle
f_*(T_{s_0}S_a)\rangle \in C_x \Xi$, if and only if $ \langle
\partial_D (-a_{1},a_{1},\ldots,-a_{k},a_{k})\rangle \in C_x \Xi
$. \QED

\remskip
Due to Proposition \ref{length}, we can restate Proposition
\ref{teoprymim} in the following way: 

\begin{cor}\label{cor615} With the same hypothesis as the proposition,
\begin{itemize}

\item[\textnormal{(a)}]
If $\MUL_x \Xi =h^0(L)/2$, then 
there is a bijection of sets
\[
\{a \in \PROJ^{k-1}\ |\  d_1(a)>0\}
\leftrightarrow
\langle\Psi_{\eta}(q_{1}),\ldots,\Psi_{\eta}(q_{k})\rangle \cap C_x \Xi.
\]

\item[\textnormal{(b)}]
If $\MUL_x \Xi >h^0(L)/2$, and 
consequently, $L=\pi^*M\otimes \OO_{\CT}(B)$, with 
$h^0(M)>h^0(L)/2$, $B\ge 0$, and $B\cap \tau^*B=\emptyset$, then
there is a bijection of sets
\[
\{a \in \PROJ^{k-1}\ |\  d_1(a)>2h^0(M)-h^0(L) \textnormal{ or } 
d_2(a) >0\}
\]
\[
\leftrightarrow
\langle\Psi_{\eta}(q_{1}),\ldots,\Psi_{\eta}(q_{k})\rangle \cap C_x \Xi.
\]
\QED
\end{itemize}
\end{cor}

We have the following elementary consequence.

\begin{cor}
With the same hypothesis as the proposition, if $\deg(D)=2k \le
2g-2$, and $h^0(L(-2D)) \ne 0$, then $C_x\Xi$ contains a
$(k-1)$-dimensional hyperplane.\QED 
\end{cor}

\begin{rem} In the above analysis in the case that both $d_1>0$
and $d_2>0$, we did not rule out the possibility that sections
lift to arbitrary order. Consequently, we did not show that
$f(S)\nsubseteq \THT$.  Nevertheless, in the case that sections
lift to arbitrary order, and hence $f(S)\subseteq \THT$, it is
clear that $f_* (T_{s_0}S) \subseteq C_x\THT$, and hence the
conclusions of the corollaries hold in these cases as well.
\end{rem} 
\remskip

We will now do a computation to prove the following 
proposition.  This will illustrate the basic technique
to be used in the next section.

\begin{pro}
Suppose that $x$ is a singular point of $\Xi$, corresponding to a
line bundle $L\in \PIC^{2g-2}(\CT)$, such that
$\MUL_x\Xi=h^0(L)=2$.  For a point $q\in C$, let
$\pi^{-1}(q)=\{p,\tt(p)\}$. Then $\Psi_{\eta}(q)\in C_x \Xi$ if
and only if $h^0(L(-\tt(p)))\ne 1$ or $h^0(L(-2p-\tt(p)))\ne 0$.
\end{pro}

\begin{rem} Since the condition $\Psi_{\eta}(q)\in C_x \Xi$ is
independent of the choice of $p$ versus $\tt(p)$, one can conclude
from the proposition that for any point $p\in C$,
$h^0(L(-\tt(p)))= 1$ and $h^0(L(-2p-\tt(p))) =0$ if and only if
$h^0(L(-p))=1$ and $h^0(L(-p-2\tt(p)))= 0$. 
\end{rem} 
\remskip

\PRF.  By Theorem \ref{teosv}, $L=\pi^*M \otimes \OO_{\CT}(B)$,
with $h^0(C,M) = h^0(\CT,L)=2$, $B \ge0$, and $B \cap \tau^*B =
\emptyset$.  Let $D=\pi^*(q)$, consider the deformation
$\LL_{D;1}$, and let $f:S\ra P$ be the associated morphism with
$f(s_0)=x$.  By Proposition \ref{length} (b), all sections lift to
first order, and I claim that if $h^0(L(-\tt(p)))\ne 1$, or
$h^0(L(-p-\tt(p)-p))\ne 0$, then a nontrivial section must lift to
second order, so that by Corollary \ref{cor615} (b), $\Psi(q)\in
C_x\Xi$. 

Indeed, assume $h^0(L(-\tt(p)))\ne 1$.  Then either $h^0(L(-p
-\tt(p)))= 2$, in which case $h^0(L(-2D))\ne 0$, and a section
lifts to second order by Corollary \ref{cor149}, or $h^0(L(-p
-\tt(p)))= 1$.  In this case, consider a nonzero section $s\in
H^0(L(-p-\tt(p)))$. Let $s +\ss^{(1)}t$ be the standard lifting of
$s$, so that the general lifting of $s$ will be of the form $s
+(\ss^{(1)}+\ffi) t$ for some $\ffi \in H^0(L)$. Then $A_2(s
+(\ss^{(1)}+\ffi) t)= (0, -(s/z) (p)-\ffi(p), 0, \ffi(\tt(p)))$.
Now considering the fact that $h^0(L(-\tt(p)))=2$ and
$h^0(L(-p-\tt(p)))=1$, it follows that the map $H^0(L(-\tt(p)))\ra
\CX$ given by $\psi \mapsto \psi(p)$ is surjective, so that there
is a section $\ffi\in H^0(L)$ such that $\ffi(\tt(p))=0$, and
$\ffi(p)=-(s/z)(p)$.  Thus, $A_2(s +(\ss^{(1)}+\ffi) t)=0$, and so
$s$ lifts to second order. 

On the other hand, suppose $h^0(L(-p-\tt(p)-p))\ne 0$, and let $s
\in H^0(L(-p-\tt(p)-p))$ be a nonzero section.  Then $s$ lifts to
first order since $s\in H^0(L(-D))$, and consequently set
$s+\ss^{(1)}t$ to be the standard first order lift. 
$A_2(s+\ss^{(1)}t) =(0, -(s/z) (p), 0, 0)$, and since $s \in
H^0(L(-p-\tt(p)-p))$ it follows that $A_2(s+\ss^{(1)}t)=0$. 

Conversely, suppose $h^0(L(-2p-\tt(p)))=0$ and
$h^0(L(-\tt(p)))=1$.  In this case I claim that only the trivial
section lifts to second order, and hence, by Corollary
\ref{cor615} (b), $\Psi_\eta(q)\notin C_x \Xi$. In fact, this
follows from the proof of Theorem \ref{teo2}; the key observation
is that the proof of the theorem depends only on the numerology of
Lemma \ref{linsys}, not on the assumption that the chosen points
were general.  Thus we must check that the conditions
$h^0(L(-2p-\tt(p)))=0$, and $h^0(L(-\tt(p)))=1$, are sufficient to
establish the results of Lemma \ref{linsys}. Using Riemann-Roch,
we need only check (b), (d), (f), (g), (h), and (j). 

(b) $h^0(L(-D))=1$; $1 \le h^0(M(-q))\le   h^0(L(-D))\le h^0(L(-\tt(p)))=1$.
(d) $h^0(M(-2q))= h^0(L(-2D))\le h^0(L(-2p-\tt(p)))=0$.
(f) $h^0(L(D+\tt(p)))=3$; this follows from Riemann-Roch. 
(g) is the same as (f) in this case. 
(h) $h^0(L(-\tt(p)))=1$ is given.
(j) is vacuous. 
\QED

\subsection{Secant varieties}\label{secsec}

We now direct our attention to secant varieties of the Prym
canonical image.  For $0\le k \le r$, the $k$-secant variety of a
curve $\Gamma$ embedded in $\PROJ^r$, is defined to be the closure of the union of the linear subspaces in $\PROJ^r$ spanned by a $(k+1)$-tuple of distinct points of $\Gamma$;
i.e. the $0$-secant variety is $\Gamma$, and the
$1$-secant variety is the usual secant variety. 

\begin{teo}\label{teosec}
Suppose that $x$ is a singular point of $\Xi$, 
corresponding to a line bundle $L\in \PIC^{2g-2}(\CT)$ such that  
$h^0(L)=2n$.
\begin{itemize}
\item[\textnormal{(a)}] 
 The $(n-1)$-secant variety of 
$\Psi_{\eta}(C)$ is not contained in 
$C_x \Xi$.
More precisely,
if 
$q_1, \ldots ,q_n$ are general points of $C$, then 
$\langle\Psi_{\eta} (q_1), \ldots ,\Psi_{\eta}(q_n) \rangle \nsubseteq
C_x \Xi$.  

\item[\textnormal{(b)}] The $(n-2)$-secant variety of 
$\Psi_{\eta}(C)$ is contained in 
$C_x \Xi$.  Hence  the 
$k$-secant variety of $\Psi_{\eta}(C)$ is contained in $C_x\Xi$ for
all $0\le k\le n-2$.

\end{itemize}

\end{teo}

\PRF. (a) Let $q_1,\ldots,q_n$ be $n$ general points of $C$, let
$D'=\sum_{i=1}^nq_i$, and let $D=\pi^*D'$.  For $a \in \CX^n$,
consider the deformation $\LL_{D;a}$, and let $f:S\ra P$ be the
associated morphism, with $f(s_0)=x$.  In the case $\MUL_x \Xi=n$,
the proof of Theorem \ref{teosv} implies that for a general $a\in
\CX^n$, $(\MUL_{s_0}\Theta_S)/2=\MUL_x \Xi$. In the case $\MUL_x
\Xi >n$, the proof of Theorem \ref{teo2} implies the same result.
Proposition \ref{teoprymim} then implies that $\langle\Psi_{\eta}
(q_{1}),\ldots,\Psi_{\eta}(q_{n})\rangle \nsubseteq C_x \Xi$.
\remskip

(b) In the case $\MUL_x \Xi=n$, for general points
$q_1,\ldots,q_{n-1} \in C$, let $D'=\sum_{i=1}^{n-1}q_i$, let
$D=\pi^*D'$, and consider the deformation $\LL_{D;a}$. For all
$a\in \CX^{n-1}$, $d_1(a) \ge h^0(L(-D)) >0$, and hence by
Corollary \ref{cor615} (a), the $(n-2)$-secant variety of
$\Psi(C)$ is contained in $C_x\Xi$. \remskip

In the case $\MUL_x \Xi>n$, if $L=\pi^*M \otimes \OO_{\CT}(B)$,
$h^0(C,M) > h^0(\CT,L)/2$, $B \ge0$, and $B \cap \tau^*B =
\emptyset$, let $n_1=h^0(M)$, and $n_2=h^0(L)-h^0(M)$. For general
points $q_1,\ldots,q_{n-1} \in C$, let $D'=\sum_{i=1}^{n-1}q_i$,
$D=\pi^*D'$, and consider the deformation $\LL_{D;a}$. We will
find that for general $a\in \CX^{n-1}$, there is a nontrivial
section which lifts to second order, and hence by Corollary
\ref{cor615} (b), the $(n-2)$-secant variety of $\Psi(C)$ is
contained in $C_x\Xi$. 

Let $a\in \CX$ be such that $a_i\ne 0$ for all $i$.  Let $s\in
H^0(L(-D))$, and let $ s + \ss^{(1)}t$ be the standard lift of
$s$, as in Lemma \ref{canon}.  Then the general lift of $s$ will
be of the form $s + (\ss^{(1)}+\ffi)t$, for some $\ffi\in H^0(L)$.
We have seen that
\[
A_2(s + (\ss^{(1)}+\ffi)t)=
(0,-a_1^2(s/z)(p_1)-a_1\ffi(p_1), 0, a_1\ffi(\tt(p_1)),\ldots).
\]

I claim that there is some $s \in H^0(L(-D))$, and some $\ffi\in
H^0(L)$ such that $A_2(s + (\ss^{(1)}+\ffi)t)=0$. From the
equation above, this is equivalent to the claim that there exists
some $s \in H^0(L(-D))$, and some $\ffi\in H^0(L)$ such that
$\ffi(p_i)=a_i(s/z)(p_i)$ and $\ffi(\tt(p_i))=0$, for $1\le i \le
n-1$. 

Let $V_1=
\{\big(a_1(s/z)(p_1),\ldots,a_{n-1}(s/z)(p_{n-1})\big)\in
\CX^{n-1} |s \in H^0(L(-D))\}$, and let $V_2=
\{\big(\ffi(p_1),\ldots,\ffi(p_{n-1})\big) \in \CX^{n-1} |\ffi \in
H^0(L (-\sum_{i=1}^{n-1}p_i))\} $.  In order to prove our claim,
we need only show that $\dim(V_1)+\dim(V_2) >n-1$. Now
$\dim(V_1)=h^0(L(-D))-h^0(L(-D-\sum_{i=1}^{n-1} p_i))$, and
$\dim(V_2)=h^0(L(-\sum_{i=1}^{n-1} \tt(p_i)))- h^0(L(-D))$. Hence,
$$\dim(V_1)+\dim(V_2)= h^0(L(-\sum_{i=1}^{n-1} \tt(p_i)))-
h^0(L(-D-\sum_{i=1}^{n-1} p_i)).$$ Since the points $p_i$ are
general, we have that $h^0(L(-\sum_{i=1}^{n-1}
\tt(p_i)))=2n-(n-1)=n+1$. Furthermore, it follows from Corollary
\ref{cor313} that $h^0(L(-D-\sum_{i=1}^{n-1} p_i))
=\max(0,n_1-n+1-(n-1))= \max(0,n_1-2n+2)\le 2$, with equality
holding if and only if $n_1=2n$.  Thus $$\dim(V_1)+\dim(V_2) \ge
n+1 -2, $$ with equality holding if and only if $n_1=2n$. 
Therefore, we have established the claim in the case $n_1 \ne 2n$,
and it follows that in this case there is a nontrivial section
which lifts to second order. 

On the other hand, the case $n_1=2n$ is much easier.  Indeed,
$h^0(L(-2D))$ $=2n - (n-1)-(n-1)=2\ne 0$, and thus as observed in
the corollary to Lemma \ref{Dn1n22}, a nontrivial section must
lift to second order. \QED

\remskip

\begin{rem} This theorem generalizes Smith-Varley \cite{sv1}
Proposition 5.1, which does not address the issue of the secant
variety, and makes the additional assumption that either
$\MUL_x\Xi=(1/2)h^0(L)$, or $L$ is base point free and
$(1/2)h^0(L)\le \MUL_x\Xi\le h^0(L)-1$.  Also, as a consequence of
Theorem \ref{teosec}, we see that we could not have proven Theorem
\ref{teo2} using a divisor $D$ of degree less than $h^0(L)$.
\end{rem}

\remskip
\begin{rem} Using similar techniques, one can easily prove the
Riemann singularity theorem for Jacobians, as well as the fact
that for a Jacobian $(JC,\Theta)$, and a point $x\in \sing\Theta$
corresponding to a line bundle $L\in \PIC^{g-1}(C)$ with
$h^0(L)=n$, the $k$-secant variety of the canonical image of the
curve is contained in $C_x\Theta$ if and only if $k\le n-2$. 
(c.f. \cite{cmf} Theorem 1.9, and \cite{acgh} Theorem 1.6, p.232.)
\end{rem}

\begin{cor}\label{cortc}
$\Psi_{\eta}(C)\subseteq C_x\Xi$ if and only if $h^0(L_x)\ge 4$.
\QED 
\end{cor}

\begin{cor}[Tjurin \cite{t}, Smith-Varley \cite{sv1}]\label{cornondeg}
If $\MUL_x \Xi =2$, then one of the following must hold:

\begin{itemize}

\item[\textnormal{(a)}]$h^0(L)=4$, and $L$ can not be written in the
form $L=\pi^*(M)\otimes \OO_{\CT}(B)$ where $h^0(M)> 2$, $B>0$ and
$B\cap \tt^*B=\emptyset$.  In this case, 
$\Psi_{\eta}(C)\subseteq C_x \Xi$;  as a result, $C_x \Xi$ 
is nondegenerate,
and $\RK(C_x \Xi)\ge 3$.  In addition, the 
secant variety of $\Psi_{\eta}(C)$ is not contained in 
$C_x\Xi$.

\item[\textnormal{(b)}]$h^0(L)=2$, and 
$L=\pi^*(M) \otimes \OO_{\CT}(B)$ where 
$h^0(M)=2$, $B\ge 0$, and $B\cap \tt^*B=\emptyset$.
In this case 
$\Psi_{\eta}(C) \nsubseteq C_x \Xi$.\QED
\end{itemize}
\end{cor}

\begin{rem}\label{remtjurin}
The fact that in (a), $\Psi_{\eta}(C) \subseteq C_x \Xi$,
was first shown by Tjurin (\cite{t} Lemma 2.3, p.963).
The fact that in (a), the secant variety of
$\Psi_{\eta}(C)$ is not contained in $C_x\Xi$, and that in (b)
$\Psi_{\eta}(C) \nsubseteq C_x \Xi$, is a consequence of Theorem
\ref{teoprymim}, and was not previously known in general.
In the special case that $(P,\Xi)$ is the Jacobian of a non hyperelliptic curve,
Smith and Varley have observed that at a generic exceptional double point $\Psi_{\eta}(C) \nsubseteq C_x \Xi$ (\cite{sv3}, p. 241, line 8).  Their argument 
in the case that the
curve has no $g^1_2$, $g^1_3$, or a $g^2_5$,
is that the Prym canonical curve is contained in every stable quadric, while by Green's theorem \cite{green} the base locus of the quadrics is a canonically embedded curve.  This can not also be a Prym canonically embedded curve, and therefore not all of the exceptional quadrics contain $\Psi_{\eta}(C)$.
\end{rem}

\subsection{Equations for Tangent Cones}

Kempf's theorem \cite{kempf} gives an equation defining the
tangent cone to $\widetilde{\Theta}$ at a point $x$ as a subscheme
of $H^0(\CT, \omega_\CT)=T_x J\CT$. In the case $T_x P \nsubseteq
C_x \widetilde{\Theta}$, this equation restricts to
$H^0(\CT,\omega_\CT)^-=T_x P$ to give the square of an equation
defining the tangent cone to $\Xi$ at $x$ as a subscheme of
$H^0(\CT,\omega_\CT)^-$. 

The aim of this section will be to give another description of the
equation for $C_x \Xi$ in the case $T_x P \nsubseteq C_x
\widetilde{\Theta}$, and to give a description of the equation of
the tangent cone in the case that $T_x P \subseteq C_x
\widetilde{\Theta}$ and $h^0(L)=h^0(M)$.  In the former case, we
will need only to look at first order liftings.  In the latter
case, we will need to look at second order liftings, but will have
the advantage of knowing that the space of sections lifting to
first order is fixed.  The analysis that follows will apply to any
situation where this is true. 

\remskip
Given $g-1$ points $q_1,\ldots,q_{g-1}$ of $C$, which are linearly
independent as points of $\Psi(C)$, let
$D=\sum_{i=1}^{g-1}\pi^{-1}(q_i)$, and consider the deformation
$\LL_{D;a}$.  Let $E_1:H^0(L(D))\ra H^0(L(D)\otimes \OO_{D})$ be
the map induced from the short exact sequence, and similarly, let
$E_2:H^0(L(2D)) \ra H^0(L(2D)\otimes \OO_{2D})$.  Let $M_{E_1}$
and $M_{E_2}$ be the matrices whose rows span the respective
images of these maps.  Define $B_1:H^0(L) \ra H^0(L(D)\otimes
\OO_{D})$ to be $A_1$. By definition, a section $s\in H^0(L)$
lifts to first order if and only if $\partial_{L;D}\circ B_1(s)=0
\in H^1(L)$. 

We will now define a map $B_2:H^0(L)\ra H^0(L(2D)\otimes
\OO_{2D})$ which will have the property that there is a one to one
correspondence between sections of $L$ which lift to second order,
and sections of $L$ in the kernel of $\partial_{L;2D}\circ B_2$.
To do this, choose a basis $\{s_1,\ldots,s_{d_1}\}$ for $W_1$, and
a set of sections $\{s_{d_1+1},\ldots,s_{2n}\}$ whose images form
a basis for $H^0(L)/W_1$.  Let
$\{s_1+\ss_1^{(1)}t,\ldots,s_{d_1}+\ss_{d_1}^{(1)}t\} \subseteq
H^0(\LL_1)$ be a set of liftings of the basis for $W_1$, and
define the map $B_2$ by $\sum_{i=1}^{2n}\aa_is_i\mapsto
A_2(\sum_{i=1}^{d_1}\aa_i(s_i+\ss_i^{(1)}t)+\sum_{j=d_1+1}^{2n}\aa_js_jt)$, where the $\aa_i\in \CX$.
One can easily check that there is a one to one correspondence
between sections of $L$ which lift to second order, and sections
of $L$ in the kernel of $\partial_{L;2D}\circ B_2$. 

Let $M_{B_1}$ and $M_{B_2}$ be matrixes whose rows span the
respective image of these maps.  Define the matrices
$M_1=\left( \begin{array}{c} 
M_{E_1} \\ 
M_{B_1} 
\end{array} \right)$
and $M_2=\left( \begin{array}{c} 
M_{E_2} \\ 
M_{B_2} 
\end{array} \right)$.

\begin{teo}
Suppose $x \in \sing \Xi$ corresponds to the line bundle
$L\in \PIC^{2g-2}(\CT)$.
\begin{itemize}
\item[\textnormal{(a)}] If $\MUL_x\Xi=h^0(L)/2$, then 
${\det (M_1)}$ is 
a homogeneous polynomial of degree $h^0(L)$, which
defines $C_x \Xi$ as a subset of $H^0(\CT,\omega_\CT)^-$.

\item[\textnormal{(b)}] If $\MUL_x\Xi=h^0(L)$, and 
$L=\pi^*M\otimes \OO_{\CT}(B)$ where $h^0(M)=h^0(L)$, 
$B\ge 0$ and $B\cap \tt^*B=\emptyset$, then 
$\det (M_2)$
is a homogeneous 
polynomial of degree $2 \cdot h^0(L)$, which defines
$C_x \Xi$ as a subset of $H^0(\CT, \omega_\CT)^-$.

\end{itemize} 
\end{teo}

\PRF.  (a) Let $\LL$ be a family of deformations parameterized by
$a\in \CX^{g-1}$, with fiber $\LL_a=\LL_{D;a}$.  By Riemann-Roch,
$h^0(L(D))=h^0(\omega_\CT \otimes L^{-1}(-D))+\deg(D)$.  But $\deg((\omega_\CT
\otimes L^{-1}(-D)))=-g-1$, and so $h^0(L(D))=2g-2$.  Letting
$\{s_1,\ldots,s_n\}$ be a basis for $H^0(L)$, we can take
$(M_{B_1})_{i;2j-1} = -a_js_i(p_j)$, and $(M_{B_1})_{i;2j} =
a_js_i(\tt(p_j))$ for $1 \le i \le n$, and $1\le j \le g-1$. 
Letting $\{r_1,\ldots,r_{2g-2-n} \}$ be a basis for
$H^0(L(D))/H^0(L)$, then we can take $(M_{E_1})_{i;2j-1}$
$=r_i(p_j)$, and $(M_{E_1})_{i;2j}=r_i(\tt(p_j))$ for $1\le i \le
2g-2-n$, and $1\le j \le g-1$.  With this notation, nontrivial
sections lift to first order if and only if $\det (M_1)=0$.  From
the form of the matrix $M_1$, it is clear that the determinant is
a homogeneous polynomial of degree $h^0(L)$ in the $a_i$.  (a) now
follows from Corollary \ref{cor615} (a). 

(b) Let $\LL$ be a family of deformations parameterized by $a\in
\CX^{g-1}$, with fiber $\LL_a=\LL_{D;a}$. In this case we have
seen that all sections lift to first order in all directions,
and it is easy to check that if $s_i+\ss_{i;j}$ is a lift
of $s_i$ in the direction of the $j$-th basis vector of 
$\CX^{g-1}$, then 
$s_i+\sum_j a_j\ss_{i;j}$ is a lift of $s_i$ in the direction of $a$.
 Using the basis
described in the definition of $B_2$, we get that
$(M_{B_2})_{i;4j-3}=0$,$(M_{B_2})_{i;4j-2}=-a_j\sum_k a_k\ss_{i;k}(p_j)$,
$(M_{B_2})_{i;4j-1}=a_j^2s_i(\tt(p_j))$ and
$(M_{B_2})_{i;4j}=a_j\sum_k a_k\ss_{i;k}(\tt(p_j))+a_j^2\frac{ds_i}{dz}(\tt(p_j))$.
Again, $M_{E_2}$ is independent of the $a_i$, and so we see that
$\det (M_2)$ is a homogeneous polynomial of degree $2h^0(L)$ in
the $a_i$.  Nontrivial sections lift to second order if and only
if $\det (M_2)=0$, and so (b) now follows from Corollary
\ref{cor615} (b). \QED

\remskip
\begin{rem} This analysis will go through in any case where the
space of sections lifting to first order is fixed.  In the case
that it is not fixed, the dependence of the entries of $M_{B_2}$
on the $a_i$ is more difficult to ascertain. 
\end{rem}

\remskip
\begin{rem} Theorem \ref{teosv} implies that in case (a), $T_xP
\nsubseteq C_x\THT$, and hence in this case Kempf's theorem gives 
an equation for the tangent cone as a scheme. 
\end{rem}

\begin{cor}[Quadric Tangent Cones]
Suppose $x \in \sing \Xi$ corresponds to the line bundle $L\in
\PIC^{2g-2}(\CT)$. If $\MUL_x \Xi =2$, then one of the following
must hold: 

\begin{itemize}

\item[\textnormal{(a)}]$h^0(L)=4$.  In this case,
$\det(M_1)=q(a_1,\ldots,a_g)^2$ for some irreducible homogeneous
quadratic polynomial $q\in \CX[a_1, \ldots,a_{g}]$. Hence,
$q=\sqrt{\det(M_1)}$ defines $C_x \Xi$ as a subscheme of
$\PROJ^{g-1}$. 

\item[\textnormal{(b)}]$h^0(L)=2$.  In this case, $\det(M_2)=q^2$
or $\ell_1\ell_2^3$, where $q$, $\ell_1$, $\ell_2$ are homogeneous
polynomials of degree two and one respectively. Hence, either $q=
\sqrt{\det (M_2)}$, or $\ell_1\ell_2$, defines $C_x \Xi$ as a
subscheme of $\PROJ^{g-1}$. 

\end{itemize}
\end{cor}

\PRF.  We have seen in Corollary \ref{cornondeg} that in case (a) the tangent cone is
nondegenerate.\QED

\remskip
\begin{rem} Smith and Varley in \cite{sv3} have used Kempf's
theorem to analyze the rank of quadric tangent cones in case (a). 
It is reasonable to expect that the description of the tangent
cone given in the corollary above will yield new information.
This is work in progress.  
\end{rem}

\begin{flushleft}
Department of Mathematics

Columbia University

New York, NY 10027

USA
\remskip
\texttt{casa@math.columbia.edu}
\end{flushleft}

\end{document}